\documentclass[12pt,leqno]{amsart}

\usepackage{amsmath,amssymb,latexsym,psfrag}
\usepackage{graphicx} 
\usepackage{amsthm,amscd,amsfonts}

\begin{document}

%\newtheorem{thm}{Theorem}
%\newtheorem{lem}[thm]{Lemma}
%\newtheorem{prop}[thm]{Proposition}
%\newtheorem{cor}[thm]{Corollary}
%\newtheorem{ex}[thm]{Example}
%\newtheorem{ddef}[thm]{Definition}
%\newtheorem{rem}[thm]{Remark}
%\newtheorem{conj}[thm]{Conjecture}

%\newtheorem{Thm}{Theorem}%[section]
%\newtheorem{Lem}{Lemma}%[section]
%\newtheorem{Cor}{Corollary}%[section]
%\newtheorem{Prop}{Proposition}%[section]
    
%% Theorem numbering and title placement:

%\newtheorem{thm}[subsection]{Theorem}
%\newtheorem{cor}[subsection]{Corollary}
%\newtheorem{lem}[subsection]{Lemma}
%\newtheorem{prop}[subsection]{Proposition}
%\newtheorem{ddef}[subsection]{Definition}
%\newtheorem{ex}[subsection]{Example}
%\newtheorem{rem}[subsection]{Remark}
%\newtheorem{conj}[subsection]{Conjecture}  

%%%%%%%% theorems etc  %%%%%%%%

\newtheorem{thm}{Theorem}[section]
\newtheorem{prop}[thm]{Proposition}
\newtheorem{lem}[thm]{Lemma}
\newtheorem{cor}[thm]{Corollary}
\newtheorem{conj}[thm]{Conjecture}
\newtheorem{ddef}[thm]{Definition}
\newtheorem{ex}[thm]{Example}
\newtheorem{rem}[thm]{Remark}
\newtheorem{notation}[thm]{Notation}

\numberwithin{equation}{section}

\newcommand{\bthm}{\begin{thm}}
\newcommand{\ethm}{\end{thm}}
\newcommand{\blem}{\begin{lem}}
\newcommand{\elem}{\end{lem}}
\newcommand{\bcor}{\begin{cor}}
\newcommand{\ecor}{\end{cor}}
\newcommand{\bprop}{\begin{prop}}
\newcommand{\eprop}{\end{prop}}
\newcommand{\bproof}{\begin{proof}}
\newcommand{\eproof}{\end{proof}}
\newcommand{\bddef}{\begin{Definition}}
\newcommand{\eddef}{\end{Definition}}

\newcommand{\beq}{\begin{equation}} \newcommand{\eeq}{\end{equation}}
\newcommand{\beqs}{\begin{equation*}} \newcommand{\eeqs}{\end{equation*}}
\newcommand{\beqarr}{\begin{eqnarray}} \newcommand{\eeqarr}{\end{eqnarray}}
\newcommand{\beqarrs}{\begin{eqnarray*}} \newcommand{\eeqarrs}{\end{eqnarray*}}
\newcommand{\barr}{\begin{array}} \newcommand{\earr}{\end{array}}
\newcommand{\btab}{\begin{tabular}} \newcommand{\etab}{\end{tabular}}

\newcommand{\bit}{\begin{itemize}} \newcommand{\eit}{\end{itemize}}
\newcommand{\ben}{\begin{enumerate}} \newcommand{\een}{\end{enumerate}}
\newcommand{\bce}{\begin{center}} \newcommand{\ece}{\end{center}}

\newcommand{\defeq}{\stackrel{\rm def}{=}}
\newcommand{\bd}{\partial}
\newcommand{\op}{\mathrm{op}}
\newcommand{\hz}{\widehat{0}}
\newcommand{\ho}{\widehat{1}}

\newcommand{\wh}{\widehat}
\newcommand{\wt}{\widetilde}
\newcommand{\cov}{\lhd}
\newcommand{\ccov}{\rhd}

%%some of the usual math blackboard caps

\newcommand{\N}{\mathbb{N}}
\newcommand{\Z}{\mathbb{Z}}
\newcommand{\C}{\mathbb{C}}
\newcommand{\R}{\mathbb{R}}
\renewcommand{\P}{\mathbb{P}}
\newcommand{\Q}{\mathbb{Q}}
\newcommand{\F}{\mathbb{F}}

%% some script caps

\renewcommand{\AA}{\mathcal{A}}
\newcommand{\BB}{\mathcal{B}}
\newcommand{\CC}{\mathcal{C}}
\newcommand{\DD}{\mathcal{D}}
\newcommand{\EE}{\mathcal{E}}
\newcommand{\FF}{\mathcal{F}}
\newcommand{\II}{\mathcal{I}}
\newcommand{\LL}{\mathcal{L}}
\newcommand{\QQ}{\mathcal{Q}}
\renewcommand{\SS}{\mathcal{S}}

%% Greek

\newcommand{\De}{\Delta}
\newcommand{\Ga}{\Gamma}
\newcommand{\de}{\delta}
\newcommand{\ga}{\gamma}
\newcommand{\la}{\lambda}
\newcommand{\La}{\Lambda}

\

\newcommand{\X}{\bf{X}}
\newcommand{\norm}[1]{\lVert#1\rVert}
\newcommand{\abs}[1]{\lvert#1\rvert}
\newcommand{\qbinom}[2]{\left[\ba{c}{#1}\\{#2}\ea\right]} %q-binom \qbinom{1}{2}
\newcommand{\comp}{\models}

\newcommand{\vanish}[1]{}

\newcommand{\bex}{\begin{ex}} \newcommand{\eex}{\end{ex}}
\newcommand{\brem}{\begin{rem}} \newcommand{\erem}{\end{rem}}

\newcommand{\bca}{\begin{cases}} \newcommand{\eca}{\end{cases}}

\def\al{\alpha}
\def\be{\beta}
\def\de{\delta}
\def\De{\Delta}
\def\ga{\gamma}
\def\Ga{\Gamma}
\def\la{\lambda}
\def\La{\Lambda}
\def\om{\omega}
\def\Om{\Omega}
\def\ze{\zeta}
\def\vphi{\varphi}
\def\veps{\varepsilon}
\def\si{\sigma}
\def\Si{\Sigma}

\def\N{\mathbb{N}}
\def\P{\mathbb{P}}
\def\Z{\mathbb{Z}}
\def\Q{\mathbb{Q}}
\def\R{\mathbb{R}}
\def\C{\mathbb{C}}
\def\L{\mathbb{L}}

\def\AA{{\mathcal{A}}}
\def\BB{{\mathcal{B}}}
\def\CC{{\mathcal{C}}}
\def\FF{{\mathcal{F}}}
\def\LL{{\mathcal{L}}}
\def\NN{{\mathcal{N}}}
\def\HH{{\mathcal H}}
\def\RR{{\mathcal R}}
\def\MM{{\mathcal M}}

\def\bb{\mathbf{b}}
\def\cc{\mathbf{c}}
\def\ff{\mathbf{f}}
\def\bg{\mathbf{g}}
\def\hh{\mathbf{h}}
\def\kk{\mathbf{k}}

\newcommand{\st}{\,:\,} 
\newcommand{\sbseq}{\subseteq}
\newcommand{\spseq}{\supseteq}
\newcommand{\larr}{\leftarrow}
\newcommand{\rarr}{\rightarrow}
\newcommand{\Larr}{\Leftarrow}
\newcommand{\Rarr}{\Rightarrow}
\newcommand{\lrarr}{\leftrightarrow}
\newcommand{\Lrarr}{\Leftrightarrow}

\def\il{\int\limits}
\def\sbs{\subset}
\def\sbseq{\subseteq}
\def\wh{\widehat}
\def\wt{\widetilde}
\def\oli{\overline}
\def\uli{\underline}
\def\langle{\left<}
\def\rangle{\right>}
\def\Lraw{\Longrightarrow}
\def\lraw{\longrightarrow}
\def\Llaw{\Longleftarrow}
\def\llaw{\longleftarrow}
\def\Llraw{\Longleftrightarrow}
\def\llraw{\longleftrightarrow}

\def\wtx{\underset{\displaystyle{\widetilde{}}}{x}}
\def\wth{\underset{\displaystyle{\widetilde{}}}{h_0}}

\def\({\left(}
\def\){\right)}
\def\no={\,{\,|\!\!\!\!\!=\,\,}}
\def\wt{\widetilde}

\def\rank{\text\rm{rank}}
\def\circm{\circmega}
\def\Lv{\left\Vert}
\def\Rv{\right\Vert}
\def\lan{\langle}
\def\ran{\rangle}
\def\wh{\widehat}
\def\no={\,{\,|\!\!\!\!\!=\,\,}}
\def\sbseq{\subseteq}
\def\circli{\circverline}
\def\aff{\text{\rm{aff}}}
\def\conv{\text{\rm{conv}}}
\def\sgn{\circperatorname{sgn}}

\def\lk{\mathrm{link}}

\def\CatTop{\sf Top}
\def\id{\rm id}
\def\colimit{\circperatorname{\sf colim}}
\def\hocolim{\circperatorname{\sf hocolim}}
\def\susp{\circperatorname{susp}}
\def\sd{\circperatorname{sd}}
\def\rk{\mathrm{rk}}
\def\sbseq{\subseteq}
\def\sbs{\subset}
\def\spseq{\supseteq}
\def\sps{\supset}
\def\ssm{\smallsetminus}

\newcommand{\LA}{L_{\AA}}
\newcommand{\FA}{F_{\AA}}
\newcommand{\CA}{C_{\AA}}
\newcommand{\MA}{M_{\AA}}
\newcommand{\ep}{\varepsilon}

\newcommand{\mob}{M\"obius function}
\newcommand{\supp}{\mathrm{supp}}
\newcommand{\card}{\mathrm{card}}

\title[Random walks, arrangements, cells, greedoids, libraries]{Random walks, arrangements, 
cell complexes, greedoids, and self-organizing libraries}
\author[Anders Bj\"orner]{Anders Bj\"orner \\ \\ \\
To L\'aszl\'o Lov\'asz on his 60th birthday}
\address{Royal Institute of Technology, Department of Mathematics,
  S-100 44 Stockholm, Sweden}
\email{bjorner@math.kth.se}
%\thanks{%{\bf Acknowledgement.}}
%\subjclass{}
%\keywords{}

\begin{abstract}
The starting point is the known fact that some much-studied random walks on permutations,
such as the Tsetlin library, arise from walks on real hyperplane arrangements. This paper
explores similar walks on complex hyperplane arrangements. This is achieved by
involving certain cell complexes naturally associated with the arrangement.
In a particular case this
leads to walks on libraries with several shelves. 

We also show that interval greedoids
give rise to random walks belonging to the same general family. 
Members of this family of Markov chains,
based on certain semigroups, have the property that all eigenvalues of the transition
matrices are non-negative real and given by a simple combinatorial formula.

Background material needed for understanding
the walks is reviewed in rather great detail. 

\end{abstract}

\maketitle

{\bf Contents}
\ben
\item[1.] Introduction
\item[2.] Real hyperplane arrangements
\bit
\item[2.1] Basics
\item[2.2] The braid arrangement
\item[2.3] Cell complexes and zonotopes
\item[2.4] The permutohedron and the $k$-equal arrangements
\eit
\item[3.] Complex hyperplane arrangements
\bit
\item[3.1] Basics
\item[3.2] Cell complexes 
\item[3.3]  Complexified $\R$-arrangements
\eit
\item[4.] Random walks 
\bit
\item[4.1] Walks on semigroups 
\item[4.2] Walks on $\R$-arrangements 
\item[4.3] Walks on $\C$-arrangements 
\item[4.4]  Walks on libraries   
\item[4.5] Walks on greedoids 
\eit
\item[5.] Appendix 
\bit
\item[5.1] A generalized Zaslavsky formula 
\item[5.2] Lattice of intervals 
\item[5.3] Interval greedoids
\eit
\een

%\begin{center}   
%\psfrag{a}{\small $a$}
%\psfrag{b}{\small $b$}
%\psfrag{4}{\small $4$}
%\includegraphics[scale=.5]{_2_1_01}
%\end{center}

\section{Introduction}

The following random walk, called {\em Tsetlin's library}, 
is a classic in the theory of combinatorial Markov chains. Consider books labeled by
the integers $1, 2, \ldots, n$ standing on a shelf in some order. A book is withdrawn
according to some probability distribution $w$ and then placed at the beginning of the shelf.
Then another book is withdrawn according to $w$ and placed at the beginning
of the shelf, and so on. This Markov chain is of interest also for computer science,
where it goes under names such as
{\em dynamic file management} and {\em cache management.}

Much is known about the Tsetlin library, for instance good descriptions of its stationary distribution, 
good estimastes of the rate of convergence to stationarity, exact formulas for
the eigenvalues of its transition matrix $P_w$, and more. 
These eigenvalues are nonnegative  real and their indexing and multiplicities, as well as their value,
are given by very explicit combinatorial data.

%The eigenvalues are as follows.

The Tsetlin library is the simplest of a class of Markov chains on permutations that
can be described in terms of books on a shelf. Instead of one customer visiting
the library to borrow one book which when returned is placed at the beginning of
the shelf, we can picture several customers who each borrows several books.
When the books are returned, the books of the first borrower are placed at the beginning
of the shelf in the induced order (i.e. the order they had before being borrowed).
Then the books of the second borrower are placed in their induced order, and so on.
Finally, the remaining books that noone borrowed stand, in the induced order,
at the end of the shelf.

%The results on Tsetlin's library have been generalized to such {\em one-shelf dynamic libraries}.
%In particular, the eigenvalues and their multiplicity turn out to have surprisingly exact
%combinatorial descriptions. 

The analysis of such a ``dynamic library''  became part of a vastly more general
theory through the work of 
Bidigare, Hanlon and Rockmore \cite{BHR}, continued and expanded by
Brown and Diaconis \cite{Bro1, Bro2, BrDi, Dia}. They created an attractive  theory of
random walks on hyperplane arrangements $\AA$ in $\R^d$,
%For general real hyperplane arrangements $\AA$, 
for which the states of the Markov chain are  the regions
making up the complement of $\cup \AA$ in
$\R^d$. When specialized to the braid arrangement, whose regions are in
bijective correspondence with the permutations of $\{1, 2, \ldots, n\}$,
their theory specializes precisely to the ``self-organizing'',
or ``dynamic'', one-shelf  library that we just described.
The theory was later further generalized by Brown \cite{Bro1, Bro2} to a class of semigroups.
\smallskip

The genesis of this paper is the question: {\em what about random walks on complex
hyperplane arrangements?} It is of course not at all clear what is meant.
The complement in $\C^d$ of the union of a finite collection of hyperplanes is a
$2d$-dimensional manifold, so what determines a finite Markov chain?

The idea is to consider not the complement itself, but  rather a certain finite cell complex 
determining  the complement up to homotopy type. In addition, we need that this
complex extends to a cell complex for the whole singularity link,
since much of the probability mass is typically placed in that extension. Such complexes
were introduced by Ziegler and the author in \cite{BjZi1}. The construction and basic properties 
partly run
parallel to a similar construction in the real case, 
well-known from the theory of oriented matroids.

The complex hyperplane walks take place on such cell complexes in a manner that will
be described in Section \ref{sect4:Carr}. These cell complexes have a semigroup
structure to which the theory of Brown \cite{Bro1} applies. Thus we get results for
complex hyperplane walks analogous to those for the real case.

As mentioned, when specialized to the real braid arrangement the general 
theory of walks on {\em real} arrangements leads
%, as shown by Bidigare, Hanlon and Rockmore \cite{BHR}, 
to the {\em one-shelf dynamic library}.
%As was initially described,
%So, in the end, what is the answer to the question that we started with:
What happens when we similarly specialize random walks on 
complex arrangements to the {\em complex}
braid arrangement? 
The pleasant answer is that we are led to  Markov chains modelling
{\em dynamic libraries with several shelves}. These are self-organizing libraries where
the books are placed on different shelves according to some  classification 
(combinatorics books, geometry books, etc.), and not only the books on each
shelf but also the shelves themselves are permuted in the steps of the Markov chain.
Depending on the distribution of probability mass there are different versions.

Here is one.
Say that a customer withdraws a subset $E\subseteq[n]$ of books from the library.
The books are
replaced in the following way. Permute the shelves so that 
the ones that contain one of the books from $E$
become the top ones,  maintaining the induced order among them and among the
 remaining shelves, which are now at the bottom. Then, on each shelf move the books 
 from $E$ to the beginning of the shelf, where they are placed in the
 induced order.

The exact description is given in Section \ref{sect4:library}. These Markov chains may be
of interest also for file management applications in computer science.
\smallskip

In this paper we take a somewhat leisurely walk through the territory
leading to complex hyperplane walks, recalling  and assembling results
along the way that in the end lead to the desired conclusions. 
%about random walks on complex hyperplane arrangements.
We are not seeking the greatest generality, the aim is rather for simplicity of statements and
illuminating ideas through special cases. Some proofs that would interfere with this aspiration are
banished to an appendix. 

\medskip

Several topics touched upon in this paper relate to joint work with L\'aszl\'o
Lov\'asz. This is the case for the $k$-equal arrangements \cite{BjLo} in
Section \ref{sect2:perm} and  for the greedoids \cite{BKL} in Section \ref{sect4:greedoids}.
It is a pleasure to thank Laci for all the pleasant collaborations and interesting discussions
over many years.

Also, I am grateful to Persi Diaconis for inspiration and encouragement, and to Jakob
Jonsson for helpful remarks.

\section{Real hyperplane arrangements}

We review the basic facts about real hyperplane arrangements.
This material is described in greater detail in many places, for instance in \cite{BLSWZ} and \cite{OrTe},
to where we refer for more detailed information. Also, we adhere to the notation for posets
and lattices in \cite{EC1}.

\subsection{Basics} \label{sect2:basics}

Let $\ell_1, \dots , \ell_t$ be %a collection of 
linear forms on $\R^d$, and 
$H_i = \{x \st \ell_i (x)=0\} \subseteq \R^d$ the corresponding hyperplanes.
We call $\AA=\{H_1, \ldots, H_t\}$ a {\em real hyperplane arrangement}.
The arrangement is {\em essential} if $\cap H_i =\{0\}$,
and we usually  assume that this is the case.

The {\em complement} $M_{\AA}= \R^d \setminus \cup \AA$ consists of
a collection $C_{\AA}$ of open convex cones $R_i$ called {\em regions}.
They are the connected components of the decomposition
$M_{\AA}= \biguplus R_i$ 
into contractible pieces.
 
 \newcommand{\LAop}{\LA^{\mathrm{op}}}
 
 With $\AA$ we associate its  {\em intersection lattice}
$L_{\AA}$, consisting of all intersections of subfamilies of hyperplanes $H_i$  
ordered by set inclusion.  Each subspace belonging to $\LA$ can be
represented by the set of hyperplanes from $\AA$ whose intersection it is.
In this way the elements of $\LA$ can be viewed either as subsets of $\R^d$ or as subsets 
 of $\AA$. The latter is for simplicity encoded as subsets of $[n]$ via the
 labeling $i \lrarr H_i$.

Let $\LAop$ denote $\LA$ with the opposite partial order, so in $\LAop$ the
subspaces of $\R^d$
%, resp. the subsets of $\AA$, 
are ordered by {\em reverse inclusion}. This is a geometric lattice, whose atoms
are the hyperplanes $H_i$.

The number of regions of $\AA$ is determined by $\LA$ via its M\"obius function
%$\mu(\,\cdot\, , \,\cdot\,)$ 
in the following way.
\begin{thm}[Zaslavsky \cite{Zas}]\label{zaslavsky}
\hspace*{2mm}
$| C_{\AA} |\, =\, \sum_{x\in L_{\AA}} |\mu (x, \ho)|
$
\end{thm}

There is a useful way to encode the position of a point $x\in \R^d$ with respect to $\AA$.
Define the {\em sign vector (position vector)} $\sigma(x)=\{\sigma_1, \dots, \sigma_t\} \in \{0,+,-\}^t$ by
$$\sigma_i \defeq\bca
0, & \mbox{if $\ell_i (x)=0$} \\
+,  &\mbox{if $\ell_i (x)>  0$} \\
-,  &\mbox{if $\ell_i (x)<   0$} \eca
$$
In words, the $i$th entry $\si_i$ of the sign vector $\si(x)$ tells us whether the point $x$
is on the hyperplane $H_i$, or on its positive resp. negative side. 

Let $F_{\AA}\defeq \sigma({\R^d}) \subseteq  \{+,-,0\}^t$
and make this collection of sign vectors into a poset by componentwise ordering via
%\bigskip

\begin{center}
\psfrag{+}{\Huge $+$}
\psfrag{-}{\Huge $-$}
\psfrag{0}{\Huge $0$}
\resizebox{!}{20mm}{\includegraphics{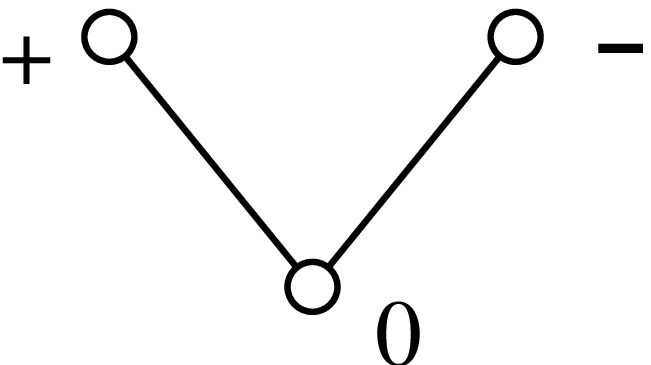}}\\
%{}{Archimedes}
\end{center}
%Note: maximal el'ts of $F_{\AA} \leftrightarrow$ regions
%\medskip

Thus, we have a surjective map $\si : \R^d \rarr \FA$.
Note that $\FA$, called the {\em face semilattice}, has minimum element $(0,\ldots ,0)$ and its
maximal elements $F_{\AA} \cap \{+,-\}^t$ are in bijective correspondence with the regions,
as is illustrated  in Figure 1.
\medskip

\begin{center} %\label{fig1}
\resizebox{!}{4.5cm}{\includegraphics{Rfaces.epsf}}
%{\bf Figure 1.} Face poset
\end{center}
\nopagebreak\vspace{.2cm} \centerline{{\bf Figure 1.} Face semilattice of an arrangement of three 
lines in $\R^2$.}
\vspace{.8cm}

The {\em composition} $X\circ Y $ of two sign vectors $X, Y \in \{0,+,-\}^t$ is
% the sign vector in $\{+,-,0\}^t$ 
defined by
%{}{Composition:} \quad $X\circ Y \in \{+,-,0\}^t$ defined by
$$ (X\circ Y)_i \defeq \bca
X_i,  \mbox{ if $X_i \neq 0$} \\
Y_i,  \mbox{ if $X_i = 0$} 
\eca
$$
This operation on $\{0,+,-\}^t$ is associative, idempotent, and has unit element $(0, \dots, 0)$.
Furthermore,  $\FA$ forms a closed subsystem: if $X, Y \in F_{\AA}$ then $ \; X\circ Y \in F_{\AA}
$.
%\beq\label{composition}
%\eeq
Here is the geometric reason: choose points $x,y\in \R^d$ such that
$\si (x)=X$ and $\si (y)=Y$. Move a small distance from $x$ along the
straight line segment from $x$ to $y$. The point $z$ reached has the
position $\si (z)=X\circ Y$.

Hence, \beq\label{semigroup}
(F_{\AA}, \circ ) \mbox{ is an idempotent semigroup.}
\eeq

The combinatorics of sign vectors  is systematically developed in {\em oriented matroid theory},
where the elements of $\FA$ are called ``covectors'' and the system $(\FA, \circ)$ is
the basis for one of the fundamental axiom systems,
see \cite[Section 3.7]{BLSWZ}.

\newcommand{\spa}{\mathrm{span}}

There is an important {\em span map}
\beq\label{spanmap}
\spa: \FA \rightarrow \LA
\eeq
which can be characterized in two ways. Combinatorially, it
sends the  sign-vector $X$ to the set of positions of its non-zero components (a subset of $[n]$).
Geometrically, it sends the cone $\si^{-1}(X)$ to its linear span.
%See Figure 2 for an illustration  based on the arrangement in Figure 1.
%\medskip

%\begin{center} %\label{fig1}
%\resizebox{!}{5cm}{\includegraphics{fig/Rfaceposet.epsf}}
%{\bf Figure 1.} Face poset
%\end{center}
%\nopagebreak\vspace{.2cm} \centerline{{\bf Figure 2.} The span map.}
%\vspace{.8cm}

The span map is  a rank-preserving and order-preserving
semigroup map, meaning that
\beqarr
\rk_{\FA} (X)&=& \rk_{\LA} (\spa(X)) \\
X \le Y &\Rarr& \spa(X)\le \spa(Y)       \\
 \spa(X\circ Y) &=& \spa(X) \vee \spa(Y) 
\eeqarr
Also, we have that
\beqarr
 X\circ Y = Y &\Lrarr& X \le Y \\
 X\circ Y= X  & \Lrarr &  \spa(Y)\le \spa(X)      
\eeqarr

\subsection{The braid arrangement} \label{sect2:braid}

\newcommand{\ZZZ}{\mathbf {Z}}
\newcommand{\perm}{\ZZZ_n^{\mathrm{perm}}}

The {\em braid arrangement} $\BB_n =\{x_i-x_j \mid 1\le i<   j\le n\}$\; in $\R^n$
plays an important role in this paper, due to its close connections with the combinatorics
of permutations and partitions.
The hyperplanes in $\BB_n$ all contain  the diagonal  line $(t,t,\ldots,t)$.
By intersecting with the hyperplane orthogonal to
this line we get an essential
arrangement, now in $\R^{d-1}$.

The intersection lattice $L_{\BB_n}$ is isomorphic to the partition lattice $ \Pi_n$,
i.e. the partitions of the set $[n]$ ordered by reverse refinement.
The correspondence between a set partition and a subspace obtained by
intersecting some hyperplanes $x_i - x_j$
is easily understood from examples:
\newcommand{\pinord}{\Pi_n^{\mathrm{ord}}}
$$(134 \mid27\mid 5 \mid 6 )\;\;\leftrightarrow\;\; \bca
x_1=x_3=x_4 \\x_2=x_7
\eca$$
and $$  \, (1345\mid 267) <      (134\mid 27\mid 5 \mid 6) .$$

\newcommand{\lang}{\langle\,}
\newcommand{\rang}{\,\rangle}
The face semilattice $F_{\BB_n}$is isomorphic to the meet-semilattice of
{\em ordered} set partitions $\Pi_n^{\mathrm{ord}}$ (so, the order of the blocks
matters), ordered by reverse refinement.
For instance,
$$ \langle\, 134\mid6\mid27\mid5 \,\rangle    \;\;\leftrightarrow\;\; \bca
x_1=x_3=x_4 \\ x_4<   x_6<   x_2 \\ x_2=x_7\\x_7<   x_5
\eca$$
and $$\lang1346 \mid257\rang  \, <    \, \lang134\mid6\mid27\mid5\rang .$$
Under this correspondence the regions of $\R^{n-1}\setminus \cup \BB_n$
are in bijection with the ordered partitions into singleton sets,
or in other words, with the permutations of the set $[n]$.
The span map (\ref{spanmap}) is the map
$\pinord\rarr\Pi_n$ that sends an ordered partition $\lang   \ldots \rang  $
 to an unordered  partition $( \ldots )$
by forgetting the ordering of its blocks.

Composition in $F_{\BB_n}$ has the following description.
If $X= \lang X_1, \dots, X_p\rang  $ and $Y=\lang Y_1, \dots, Y_q\rang  $ 
are ordered partitions of $[n]$,
%\; $X_i, Y_j \subseteq [n]$,
then $X\circ Y=\lang  X_i \cap Y_j\rang  $ with the blocks ordered by
the lexicographic order of the pairs of indices $(i,j)$.
For instance,
$$ \lang257\mid 3\mid 146\rang   \circ  \lang17\mid 25 \mid 346\rang   = \lang7\mid 25\mid
3\mid 1\mid 46\rang ,
$$
as can conveniently be seen from the computation table
\beq\label{state0}
\barr{c||c|c|c||}
\circ &1,7 & 2,5 &3,4,6 \\
\hline\hline
2,5, 7& 7 &2,5& \\
\hline
3& &  &  3\\
\hline
1,4,6& 1&&4,6 \\
\hline\hline
\earr
\eeq

\vspace{2mm}

\subsection{Cell complexes and zonotopes} \label{sect2:cells}

The whole idea of random walks on complex hyperplane arrangements rests on the idea of
walking on the cells of an associated cell complex. We therefore review the construction
used in \cite{BjZi1} of such cell complexes. The basic idea is given together with
two applications. The first one is the construction of cell complexes for the complement
of a linear subspace arrangement in $\R^d$ at the end of this section. The other is the
construction of cell complexes for hyperplane arrangements in $\C^d$, to
which we return in Section \ref{sect3:cells}. See e.g. \cite{Bjo} for topological terminology.

A regular cell decomposition $\Ga$ of the unit sphere $S^{d-1}$ is said to be  {\em PL} if its
barycentric subdivision (equivalently, the order complex of its face poset) is a
piecewise linear triangulation of $S^{d-1}$. Here is a simple combinatorial
procedure for producing regular cell complexes of 
certain specific homotopy types  from posets.

\bprop{\cite[Prop. 3.1]{BjZi1}}\label{cells1}
Suppose that $\Ga$ is a PL regular cell decomposition of $S^{d-1}$, with face poset
$F_{\Ga}$. Let $T\subseteq S^{d-1}$ be 
%a subcomplex, meaning that $T$ is 
a subspace of the sphere such that  $T=\cup_{\tau\in G} \tau$
for some order ideal $G\subseteq F_{\Ga}$. Then 
the poset $(F_{\Ga} \setminus G)^{\mathrm{op}}$ 
%with the opposite order
is the face poset of a regular cell complex having the homotopy type of the 
complement $S^{d-1}\setminus T$.
\eprop

Now, let $\AA$ be an essential hyperplane arrangement in $\R^d$.
For a general sign vector $X\in \FA$ the set $\si^{-1}(X)$ is a convex cone in $\R^d$ which is open in its
linear span. Let $\tau_X \defeq \si^{-1}(X)\cap S^{d-1}$.
The sets $\tau_X$, for $X\in \FA\setminus \hz$, partition the 
 the unit sphere and are in fact the open cells of a regular CW
decomposition of $S^{d-1}$. Furthermore,
the inclusion relation of their closures $\overline{\tau_X}$ coincides with the partial order
we have defined on $\FA$. Thus, $F_{\AA}\setminus \hz$ is the face poset of a regular cell decomposition 
$\Ga_\AA$ of the unit 
sphere in $\R^d$, namely the cell decomposition naturally cut out by the hyperplanes.

The cell complex $\Ga_\AA =\{\tau_X\}_{X\in \FA}$ induced by a hyperplane arrangement
$\AA$ is PL. Thus, via Proposition \ref{cells1} we can construct cell complexes determining the
complement of a subcomplex up to homotopy type. Combinatorially the description
is simple: erase from the face poset $\FA$ all the cells that belong to the given subcomplex
and then turn the remaining subposet upside down. Done!

\newcommand{\ZA}{\ZZZ_{\AA}}

The cell complexes constructed this way from a
hyperplane arrangement  $\AA$   can be geometrically realized
on the boundary of an associated convex polytope. Namely,
with  $\AA$  is associated its {\em zonotope}
$\ZZZ_{\AA} = [-e_1, e_1] \oplus \cdots \oplus [-e_t, e_t]$. Here $e_i$ is a normal vector in $\R^d$
to the hyperplane $H_i$ and the right-hand side denotes Minkowski sum
of centrally symmetric line segments. Thus, $\ZA$ is a centrally symmetric convex polytope,
determined this way up to combinatorial equivalence. A key property
of $\ZA$ is that  there exists an order-reversing bijection
between the  faces on its boundary and the cells of $\Ga_\AA$.
In other words,  the poset of proper faces of  $\ZA$ is isomorphic to the opposite
of the face poset of $\AA$: 
\beq\label{dual}F_{\,\ZA}\, \cong\, (\FA \setminus \hz\, )^{\op}\eeq
\smallskip

Suppose that $\AA$ is an arrangement of linear subspaces
of arbitrary dimensions in $\R^d$. Say that we want to construct a cell complex 
having the homotopy type of its complement  
$\R^d\setminus \cup\AA$. This complement is by radial projection homotopy equivalent to
its intersection with the unit sphere $S^{d-1}$. Therefore 
the preceding construction is applicable. We just have
to choose an auxiliary  hyperplane arrangement $\HH$ into which $\AA$ {\em embeds},
meaning that each subspace in $\AA$ is
the intersection of  some of the hyperplanes from $\HH$. This is clearly always possible.
Putting the various pieces of information together and applying Proposition \ref{cells1}
we obtain the following description.

\bthm{\cite{BjZi1}}\label{cells2}
Let $\AA$ be an arrangement of linear subspaces in $\R^d$.
Choose a hyperplane arrangement $\HH$ into which $\AA$ embeds.
Then the complement $\R^d\setminus \cup\AA$ has the homotopy type of a
subcomplex $\ZZZ_{\HH, \AA}$ of the boundary of the zonotope $\ZZZ_\HH$.
The complex $\ZZZ_{\HH, \AA}$ is obtained by deleting from the
boundary of $\ZZZ_\HH$ all faces that correspond to cells $\tau_X$ contained in $\cup \AA$.
\ethm

\subsection{The permutohedron and the $k$-equal arrangements}  \label{sect2:perm}

\newcommand{\ZZ}{{\bf Z}}
\renewcommand{\SS}{{\bf S}}

We illustrate the general constructions of the preceding section by applying them to
 the {\em $k$-equal arrangements} $\AA_{n,k}=
\{  x_{i_1}= x_{i_2} = \cdots = x_{i_k}\st 1\le i_1<   i_2<    \cdots <   i_k\le n    \}$ in $\R^n$.
The topology of their complements play a crucial role in the solution
of a complexity-theoretic problem in joint work with Lov\'asz and Yao \cite{BjLo, BLY}.
See also \cite{BjWe}, where their homology groups were computed.
The $k$-equal arrangements  embed into the braid arrangement (the $k=2$ case), so 
Theorem \ref{cells2} is applicable.
It tells us that, up to homotopy type, the topology of the complement of the $k$-equal
arrangement $\AA_{n,k}$ is realized by some subcomplex of the 
zonotope of the braid arrangement. This subcomplex can be very explicitly described.

\begin{center} %\label{fig1}
\resizebox{!}{6cm}{\includegraphics{perm.epsf}}
%{\bf Figure 1.} Face poset
\end{center}
\nopagebreak\vspace{.2cm} \centerline{{\bf Figure 2.} The permutohedron $\ZZZ_4^{\mathrm{perm}}$.}
\vspace{.8cm}

The zonotope of the braid arrangement  $\BB_n$ is the {\em permutohedron} $\perm$,
that is, the convex hull of the $n!$ points in $\R^n$ whose coordinates are
given by a permutation of the numbers $1, 2, \ldots, n$. Its $n!$ vertices are
in bijection with the $n!$ regions of $\BB_n$, in accordance with the duality (\ref{dual}).

%See e.g. \cite[Section 2.2]{BLSWZ} for details
%and references concerning this construction. 

 We want to describe the subcomplex
$\ZZ_{\AA_{n,k}}$ of the boundary of $\perm$ which is homotopy
equivalent to the complement $\MM_{n,k}$ of $\AA_{n,k}$.

For this one argues as follows, keeping Section \ref{sect2:braid} in fresh memory. Let
$f: \Pi_n^{\mathrm{ord}} \rightarrow \Pi_n$
be the span map, i.e., the forgetful map that sends an ordered partition of
$[n]$ to the corresponding unordered partition.
The set $\Pi_n^{\mathrm{ord}}\setminus \hz$ ordered by  refinement is the 
poset of proper faces 
of the permutohedron $\perm$, whereas
the set $\Pi_n$ ordered by refinement is the opposite of intersection
lattice of the braid arrangement. 
The image $f(\pi )$ for 
$\pi\in \Pi_n^{\mathrm{ord}}$ is a partition determining the 
span of the corresponding cell (i.e., the smallest intersection
subspace of the braid arrangement in which the cell is
contained). More precisely, the span of  $\pi $
is the subspace obtained by setting $x_{i_{1}} 
=x_{i_{2}}=\dots =x_{i_{j}}$ for each block $\{ i_{1}, 
i_{2}\dots ,i_{j} \}$ of $\pi$. Thus, a cell
$\pi\in \Pi_n^{\mathrm{ord}}$ lies in the union of the $k$-equal
arrangement if and only if  some block has size at
least $k$.

It follows that the complex $\ZZ_{\BB_n, \AA_{n,k}}$ consists of
those cells on the boundary of the permutohedron
$\perm$ that correspond to
ordered partitions with all blocks of size less than
$k$. If an ordered partition has blocks of sizes
$b_1, \dots , b_e$, then the corresponding face of $\perm$ 
is the product of smaller permutohedra of dimensions
$b_1 -1, \dots , b_e -1$. Therefore, the final description
of the cell complex $\ZZ_{\BB_n \AA_{n,k}}$ is that one should
delete from $\perm$ all faces that contain a $q$-dimensional
permutohedron, for $q \geq k-1$, in its decomposition.

We are led to the following result, obtained independently by E. Babson
for $k=3$ (see \cite{BBLL}) and the author \cite{Bjo01}.

\bthm\label{keq}
Delete from the boundary of the permutohedron $\perm$ every face
that contains a $d$-dimensional permutohedron, $d\ge k-1$,
in its decomposition.
Then the remaining subcomplex has the homotopy type of the complement of
the $k$-equal arrangement.
\ethm

Thus, for $k=2$ one  deletes everything but 
the vertices, for $k=3$ one deletes all cells
except those that are products of edges (equivalently,
keep only the cubical faces),
for $k=4$ one deletes all cells
except those that are products of edges ($1$-dimensional zonotope) 
and hexagons ($2$-dimensional zonotope), and so on.

The case $k=3$ is especially interesting. The complex is in that
case cubical. In particular, the fundamental group of $\MM_{n,3}$
is the same as the fundamental group of the cell complex obtained 
from the graph ($1$-skeleton) of $\perm$ by gluing a $2$-cell (membrane)
into every $4$-cycle.

\brem{\rm  \label{codim2}
What was just said is part of a
more general result about gluing $2$-cells into $4$-cycles of a zonotopal
graph. 

Let $\HH$ be an arbitrary central and essential
hyperplane arrangement, and let $\AA$ be the subspace arrangement
consisting of codimension $2$ intersections of $3$ 
or more planes from $\HH$ (assuming that there are such).

Next, let $G$ be the $1$-skeleton of the
zonotope $\ZZ_{\,\HH}$. The $2$-cells
of $\ZZ_{\,\HH}$ are $2m$-gons (corresponding to codimension $2$ subspaces where
$m$ planes meet). Let $\Ga_{\AA}$ be the cell complex obtained by
gluing $2$-cells into the $4$-cycles of the graph $G$.
Then the general construction
above shows (since fundamental groups live on $2$-skeleta)
that the fundamental group of $\Ga_{\AA}$
is isomorphic to that of the complement $\MM_{\AA}$.

One can go on and describe the higher-dimensional cells needed
to obtain a cell complex having the homotopy type of
the complement of such a codimension $2$  arrangement $\AA$.
They are all the cubes in the boundary of $\ZZ_{\,\HH}$, just like for
the special case of the $3$-equal arrangement.}

\erem
\brem{\rm
The two-dimensional faces of $\perm$ are either $4$-gons or $6$-gons.
What happens if we take the graph of the permutohedron and glue in
only the hexagonal $2$-cells? The answer is that we get a two-dimensional
cell complex whose 
fundamental group is isomorphic to that of the complement
of another subspace arrangement, namely the arrangement $\AA_{[2,2]}$ consisting  of
codimension $2$ subspaces of $\R^n$ obtained as intersections of pairs of hyperplanes
$x_i=x_j$ and $x_k=x_l$, for all distinct $i,j,k,l$.
Actually, for $\AA_{[2,2]}$ a stronger statement
is true:  the $2$-dimensional  cell complex
described (i.e. the permutohedron graph plus all hexagonal $2$-cells)
has the homotopy type of its complement.}
\erem

It is an interesting fact that the codimension $2$ arrangements 
%we have discussed in this section --- 
${\AA_{n,3}}$ and
$\AA_{[2,2]}$, corresponding to the two ways of gluing $2$-cells into
the permutohedron graph, share a significant
topological property, namely that their complements are K$(\pi,1)$ spaces.
See Khovanov \cite{Kho}.

\section{Complex hyperplane arrangements} 

We now move the discussion to complex space. To begin with many
of the concepts and results are parallel to the real case. But new
interesting features soon start to appear. 
This whole chapter summarizes material from \cite{BjZi1}.

\subsection{Basics} \label{sect3:basics}

%We begin with a rundown of  basic concepts and notation.

We call $\AA=\{H_1, \ldots, H_t\}$ a {\em  complex hyperplane arrangement} if
$H_i = \{z \st \ell_i (z)=0\} \subseteq \C^d$ for some 
linear forms $\ell_1, \dots , \ell_t$ on $\C^d$. 
%In what follows a particular choice of linear forms for $\AA$ is often tacitly assumed.
A particular choice of defining linear forms  %$\ell_1, \dots , \ell_t$ 
is assumed throughout, so we can also write $\AA=\{\ell_1, \ldots, \ell_t\}$.
The arrangement is {\em essential} if $\cap H_i =\{0\}$,
and we usually  assume that this is the case.
The real and imaginary parts 
of  $w=x+iy \in \C$ are denoted, respectively, by $\Re(w)=x$ and $\Im(w)=y$.

The position of a point $z\in \C^d$ with respect to $\AA$  
is combinatorially encoded in the following way.
Define the {\em sign vector (position vector)}
$\sigma(z)=\{\sigma_1, \dots, \sigma_t\} \in \{0,+,-,i,j\}^t$ by
$$\sigma_i =\bca
0,  &\mbox{if $\ell_i (z)=0$} \\
+,  &\mbox{if $\Im(\ell_i(z))=0$, $\Re(\ell_i (x)>  0$} \\
-,  &\mbox{if $\Im(\ell_i(z))=0$, $\Re(\ell_i (x)<   0$} \\
i,  &\mbox{if $\Im(\ell_i(z))>  0$} \\
j,  &\mbox{if $\Im(\ell_i(z))<   0$} 
\eca$$
%In words, the $i$th entry $\si_i$ of the sign vector $\si(z)$ tells us whether the point $z$
%is on the hyperplane $H_i$, or on its positive resp. negative side. 
Let $F_{\AA}\defeq \sigma({\C^d}) \subseteq  \{0,+,-,i,j\}^t$
and make this collection of sign vectors into a poset, called the {\em face poset},
by componentwise ordering via
\smallskip

%\\  {.} \hspace{12cm}
\begin{center}
\psfrag{i}{\Huge $i$}
\psfrag{j}{\Huge $j$}
\psfrag{+}{\Huge $+$}
\psfrag{-}{\Huge $-$}
\psfrag{0}{\Huge $0$}
\resizebox{!}{30mm}{\includegraphics{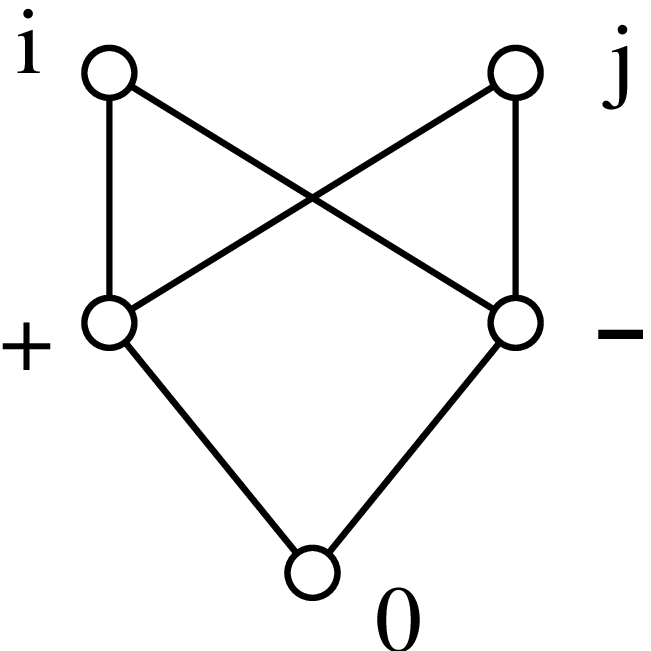}}\\
%{}{Archimedes}
\end{center}
%The face poset $\FA$ has the following combinatorial properties.
\bprop{\cite{BjZi1}}\label{posetprop}
\ben
\item $\FA$ is a ranked poset of length $2d$. Its unique minimal element is $0$.
\item The maximal elements of $\FA$ are the sign vectors in $\FA \cap \{i,j\}^{t}$.
\item $\mu(Z,W)=(-1)^{\rk(W)-\rk(Z)}$, for all $Z\le W$ in $\FA\cup \ho$.
\een
\eprop
Figure 3 (borrowed from \cite{BjZi1}),
shows the face poset of  $\AA=\{z, w, w-z\}$ in $\C^2$. The reason for
marking the elements not containing any zero
with filled dots becomes clear  in Section \ref{sect3:cells}

\begin{center} %\label{fig1}
\resizebox{!}{7cm}{\includegraphics{Cfaces.epsf}}
%{\bf Figure 1.} Face poset
\end{center}
\nopagebreak\vspace{.2cm} \centerline{{\bf Figure 3.} Face poset of an
 arrangement of three lines in $\C^2$.}
\vspace{.8cm}

The {\em composition} of two complex sign vectors  $Z\circ W \in \{0,+,-,i,j\}^t$ is defined by
\beq\label{Ccomp}
 (Z\circ W)_i = \bca
Z_i, &\mbox{if $W_i \not >    Z_i$} \\
W_i, &\mbox{if $W_i >    Z_i$} 
\eca
\eeq
Just as in the corresponding real case
this operation on $\{0,+,-,i,j\}^t$ is associative, idempotent, and has unit element $(0, \dots, 0)$.
Also, for geometric reasons (analogous to the ones in the real case)
$X, Y \in F_{\AA}$ implies that $X\circ Y \in F_{\AA}$
Hence,
\beq\label{semigroup2}
(F_{\AA}, \circ ) \mbox{ is an idempotent  semigroup.}
\eeq

For complex arrangements the notion of intersection lattice splits into two.

{\bf 1.} The  {\em intersection lattice}
$L_{\AA}$ consists of all intersections of subfamilies of hyperplanes $H_i$  
ordered by set  inclusion. 
%As in the real case, this is a geometric lattice whose atoms
%are the hyperplanes $H_i \in \AA$. 
%Each subspace belonging to $\LA$ can be
%represented by the set of hyperplanes from $\AA$ whose intersection it is.
%So, the elements of $\LA$ can be viewed either as subsets of $\C^d$ or as subsets 
 %of $\AA$. 
 %The latter can more simply be encoded as subsets of $[n]$ via the labeling $i \lrarr H_i$.

\newcommand{\Aaug}{\AA_{\mathrm{aug}}}
\newcommand{\LAaug}{L_{\AA, \,\mathrm{aug}}}
\newcommand{\LAaugop}{L_{\AA, \,\mathrm{aug}}^{\mathrm{op}}}

{\bf 2.} The {\em augmented intersection lattice} $\LAaug$ % $\LA_{\mathrm{aug}}$ of $\AA$
is the collection of all intersections of subfamilies of the {\em augmented arrangement}
$$\Aaug=\{H_1, \ldots, H_t, H_1^{\R}, \ldots, H_t^{\R}\}$$
ordered by set inclusion.
Here, $H_i^{\R}\defeq \{z\in \C^d \st \Im(\ell_i(z))=0\}$ is a
$(2d-1)$-dimensional real hyperplane in $\C^d \cong \R^{2d}$ containing $H_i$.

Again as in the real case, we denote by $\LAop$ and $\LAaugop$ the
opposite lattices, obtained by reversing the partial order.
\bprop%{\cite{BjZi1}}
\label{semimodular}
\ben
\item $\LAop$ is a geometric lattice of length $d$. 
\item $\LAaugop$ is a semimodular lattice of length $2d$. 
\een
\eprop

There is a {\em span map}
\beq\label{Cspan}\spa : \FA \rightarrow \LAaug
\eeq
defined by sending  the convex cone $\si^{-1}(Z)$, for $Z\in \FA$, to the intersection of all
subspaces in $\Aaug$ that contain $\si^{-1}(Z)$.

This map preserves poset and semigroup structure as well as poset rank.
%has these  important properties:
\bprop{\cite{BjZi1}}\label{aug}
\beqarr
\rk_{\FA} (Z)&=& \rk_{\LAaug} (\spa(Z)) \\
Z \le W &\Rarr& W \circ Z= W  \; \;\Lrarr \;\;  \spa(Z)\le \spa(W)       \\
 \spa(Z\circ W) &=& \spa(Z) \vee \spa(W) 
\eeqarr
\eprop

\subsection{Cell complexes} \label{sect3:cells}

The {\em complement} $M_{\AA}= \C^d \setminus \cup \AA$  is a complex manifold of
real dimension $2d$. There is a huge literature on the topology
of such spaces, see e.g. \cite{OrTe}.
Among the basic results we mention that
the Betti numbers of $\MA$ are determined by $\LA$ via its M\"obius function
%$\mu(\,\cdot\, , \,\cdot\,)$ 
in the following way.
\begin{thm}{\cite[p. 20]{OrTe}}\label{AOS}
\quad
$\mbox{$\beta_i (\MA) =\sum_{x\in \LA \; : \; \rk (x)=d-i} |\mu (x, \ho)|$}
$
\end{thm}

Let $\AA$ be an essential complex hyperplane arrangement in $\C^d$, as before.
For every sign vector $Z\in \FA\setminus 0$ the inverse image $\si^{-1}(Z)$ is a 
relative-open convex cone in $\C^d$. The intersections of these cones with the unit
sphere $S^{2d-1}$ in $\C^d$ are the open cells of a PL regular cell decomposition of
$S^{2d-1}$ whose face poset is isomorphic to $\FA$. Hence, as an application 
of Proposition \ref{cells1} we get part (3) of the following result. Part (2) can be
seen from the fact that  $x$ is an $\rk(x)$-dimensional linear subspace,
so $x \cap S^{2d-1}$
 is an $(\rk(x)-1)$-dimensional sphere,
for all $x\in \LAaug \setminus \hz$, where 
``$\rk$'' denotes poset rank in $\LAaug$.

\bthm{\cite{BjZi1}} \label{cells5}
\ben
\item  The poset $F_{\AA}$ is the face poset of a regular cell decomposition of the unit 
sphere in $\R^{2d}\cong \C^d$.
\item The subposet $\spa^{-1}((\LAaug)_{\le x})$  is the face poset of a regular cell decomposition of 
the sphere $S^{\rk(x)-1}$, for all $x\in \LAaug \setminus \hz$.
\item The subposet $\CA\defeq \FA \cap  \{+,-,i,j\}^t$, with opposite order, is the face poset of a regular cell
complex having the homotopy type of the complement $\MA$.
\een
\ethm

For an example, have a look at Figure 3. The sign vectors in $\FA$ that lack a zero 
component are shown by 
filled dots. Hence, the cell complex $C_{\AA}$ can be viewed
by turning the page upside-down and looking at the subposet of filled dots only.

Combining some of  this topological information with Theorem \ref{genzaslavsky} 
of the Appendix
%to the span map $ : \FA \rightarrow \LAaugop$
we obtain the following analogue of Zaslavsky's theorem  \ref{zaslavsky}
for the number of maximal cells
in the complex case.
\begin{thm}\label{genzaslavskythm}
\quad $| \max(\FA) |\, =\, \sum_{x\in \LAaug} |\mu (x, \ho)|
$
\end{thm}
\bproof
We apply Theorem \ref{genzaslavsky} 
to the span map $ : \FA \rightarrow \LAaug$. There are six conditions to verify.
%We must verify the conditions (1)--(6) of the Corollary.
With the exception of (5), they all follows from Propositions \ref{posetprop} and \ref{semimodular}.
Condition (5) is the consequence for the Euler characteristic of
%the fact that $x\in\LAaug$ is itself  a $\rk(x)$-dimensional linear subspace,
%hence its intersection with the unit sphere is a $(\rk(x)-1)$-dimensional sphere, and
%$\spa^{-1}((\LAaugop)_{\le x})$
%is the face poset of a cell complex decomposing this sphere. 
Theorem \ref{cells5}(2).

\eproof

\subsection{Complexified $\R$-arrangements}   \label{sect3:complexreal}

This section concerns the special case when all the linear forms $\ell_i(z)$ 
have real coefficients. The forms then define  both a real arrangement $\AA^{\R}$ in $\R^d$
and a complex arrangement $\AA^{\C}$ in $\C^d$. These are of course related, and we 
here summarize what expression this relation takes for the combinatorial structures
of interest.

%We begin with the fact that the  face poset $F_{\AA^{\C}}$ is determined by $F_{\AA^{\R}}$.
First a few observations about complex sign vectors. A sign vector $Z$ is called {\em real}
if all its entries come from $\{0, +, -\}$. %and {\em imaginary} if they all come  from $\{0,i,j\}$.
Every complex sign vector $Z$ can be obtained as a composition $Z=X\circ iY$ 
\footnote{\ here $i\cdot 0=0,$\; $i\cdot + = i,$\; $i\cdot - =j$.}
for two real sign vectors $X$ and $Y$. Only the vector $Y$ is unique in this decomposition.

For any poset $P$, let Int($P$) denote the set of its closed intervals.
%$[p,q]=\{r\st p\le r\le q\}$.
In the case of the face poset $F_{\AA^{\R}}$ of a real arrangement $\AA^{\R}$ we make 
$\mathrm{Int}(F_{\AA^{\R}})$ into a poset by 
introducing the following partial order:
\beq\label{Cintorder}
 [Y,X] \le [R,S] \;\;\leftrightarrow\;\; \bca
Y\le R \\ R\circ X\le S \eca
\eeq

\bprop{\cite{BjZi1}}
The map
$ \phi: \mathrm{Int}(F_{\AA^{\R}})  \rightarrow F_{\AA^{\C}}$ given by
$[Y,X]  \mapsto  X\circ iY$
is a poset isomorphism.
\eprop

For example,
$$\phi: \;  [\, (0\,- +\, 0 \,0\, -) \, ,\, (--+\, 0\,+-)\, ] \;\; \mapsto \;\;( -\, j\,\, i\,0\,+j\,)$$

Hence, the entire structure of the complex face poset $F_{\AA^{\C}}$ can be dealt with in terms
of intervals in the real face poset $F_{\AA^{\R}}$. 
In particular, the cells in the complement of $\AA$,
being the sign vectors  without any zero coordinate, 
get this description. 
$$C_{\AA^{\C}} \;\; \stackrel{\phi}{\leftrightarrow}\;\; \mbox{intervals $[Y,X]$ with } X\in\max(F_{\AA^{\R}})
$$
Composition of complex sign vectors (\ref{Ccomp}) takes the following form when translated
to intervals:
\beq\label{Ccompint}
 [Y,X] \circ [R,S]  = [Y\circ R, \; Y\circ R\circ X\circ S] 
\eeq

The augmented intersection lattice $L_{\AA^\C, \text{}aug}$ is similarly determined by
the intervals of $L_{\AA^{\R}}$, namely
\beq\label{Cintersect}
L_{\AA^\C, \text{}aug} \cong \text{Int}(L_{\AA^{\R}}),
\eeq
this time with the partial order %of intervals in $L_{\AA^{\R}}$ 
defined by
$$(x,y) \le (x', y') \quad \mbox{ if and only if \quad $x\le x'$ and }y\le y' .
$$
The span map is the natural one
\beq\label{Cspanint}
 \mathrm{Int}(F_{\AA^{\R}})  \cong F_{\AA^{\C}}  \rightarrow
 L_{\AA^\C, \text{}aug} \cong \text{Int}(L_{\AA^{\R}})
\eeq
sending $[Y,X]$ to $[\spa(Y), \spa (X)]$.
The \mob\ of $\text{Int}(L_{\AA^{\R}})$ is described in
terms of the \mob\ of the lattice $L_{\AA^{\R}}$ in
 Appendix \ref{sectApp:intervals}.
\smallskip 
%\subsection{The braid arrangement}  \label{sect3:braid}

\bex {\rm
The braid arrangement
$\BB_n^{\C}=\{x_i-x_j \mid 1\le i<   j\le n\}$\; in $\C^n$ is the complexification of the
real braid arrangement, discussed in Section \ref{sect2:braid}.
Hence we can  translate  its combinatorics into the language of intervals, as
outlined in this section.

We obtain that $\BB_n^{\C}$ has  face semilattice $$F_{\BB_n^{\C}} \cong 
\mathrm{Int}(F_{\BB_n^{\R}}) \cong\mathrm{Int}(\pinord)$$
and augmented intersection lattice
$$L_{\BB_n^{\C}, \text{aug}} \cong 
\mathrm{Int}(L_{\BB_n^{\R}}) \cong\mathrm{Int}(\Pi_n).$$

Thus, the complex sign vectors of $\BB_n^{\C}$ are encoded into pairs
$[Y,X]$ of ordered partitions, where $X$ is an refinement of $Y$.
The composition (\ref{Ccompint}) is  illustrated in this computation table:
\beq\label{state4}
\barr{c||c|c|c||c|c||c|c||}
\circ &1  &3  &5&4  & 7 & 6 & 2 \\
\hline\hline
3, 7&  & 3 &  & &7  & & \\
\hline
1&1  &  &  &  & &  &  \\
\hline\hline
2,5,6&  &  &5  &  &  & 6&2  \\
\hline
4&  &  &  & 4 & &  &  \\
\hline\hline
\earr
\eeq
from which we read that
$$ \lang 37 \mid 1 \mid\mid 256\mid 4\rang   \circ \lang 1\mid3\mid5\mid\mid4\mid7\mid\mid6\mid2
\rang  = \lang 3\mid 1\mid\mid7\mid\mid 5\mid\mid 4\mid\mid6\mid2 \rang 
$$
Here single bars denote the separation of the ground set $[7]$ into ordered 
blocks according to $X$, and double bars the coarser partition $Y$.
The rule is to read off the coarser partition of the composition by ordering
the double bar boxes lexicographically, and then read off the refinement by
ordering the single bar boxes within each
double bar box  lexicographically (empty boxes are skipped).

Notice that the cells in the complement of the complex braid arrangement, cf. Theorem \ref{cells5} (3), 
correspond to block-divided permutations:
\beqarrs
C_{\BB_n^{\C}} &\lrarr&
\mbox{sign vectors $X\circ iY$ without zero coordinates}\\
&\lrarr&
\text{intervals $[Y,X]$, $X$ maximal}\\ 
&\leftrightarrow& \mbox{permutations $X$ divided into ordered blocks $Y$}
\eeqarrs

}
\eex

\section{Random walks}
This chapter begins with a summary of Brown's theory for random walks on a 
class of semigroups \cite{Bro1}. The motivating example, namely walks on real hyperplane
arrangements, is then recalled. % in Section \ref{sect4:Rarr}.
After that  comes a sequence of applications.

\subsection{Walks on semigroups} \label{sect4:semigroups}

A {\em semigroup} is a set with an  associative composition. We also assume the existence
of an identity element,  denoted ``$e$'', and we write the composition in
multiplicative notation.

\begin{ddef}
%\bddef: 
An {\em LRB semigroup} is a finite semigroup $\Si$ with identity satisfying
\ben
\item $x^2=x  \mbox{ for all $x\in \Si$}$,
\item $xyx=xy  \mbox{ for all $x,y \in \Si$}$.
\een
A {}{\em left ideal} of $\Si$ is a subset $I\subseteq \Si$ such that $x\in \Si$, $y\in I$ $\Rightarrow$ $xy\in I$.
%\eddef
\end{ddef}

The acronym LRB stands for ``Left-Regular Band'', a name by which this class of semigroups
is sometimes known in the literature. Brown \cite{Bro1} defined 
a class of random walks on semigroups of this type.
%which is the whole purpose  for getting involved with them here. 
This section summarizes some
material from \cite{Bro1}, to where we refer for more information, background
and references.

\newcommand{\mSi}{\max(\Si)}

\begin{ddef}
Let $I$ be a left ideal of $\Si$, and let  $w$ be a probability distribution on $\Si$.
A random walk on $I$ is defined in the following way. If the current position of the
walk is at an element $y\in I$, then choose $x\in\Si$ according to the distribution $w$
and move to $xy$. \label{randomwalk}
\end{ddef}

Brown's main theorem gives surprisingly exact information about such random walks. 
In order to be able to state
it we need to first introduce two related poset structures.

Let $\Si$ be an LRB semigroup. We define a relation ``$ \,\le\, $'' on $\Si$ by
\beq x\le y \quad \Leftrightarrow \quad xy=y
\eeq
This turns out to be a partial order relation, so we may think of an LRB semigroup also as
a poset. The identity element $e$ is the unique minimal element. The set
$\mSi$ of maximal elements is a left ideal in $\Si$.

There is also another partial order significantly  related to $\Si$.

\bprop[\cite{Bro1}] \label{Broprop}
Let $\Si$ be an LRB semigroup. Then there exists a unique finite lattice $\La$ 
and an order-preserving and surjective map
\beq \supp: \, \Si \,\rarr\, \La 
\eeq
such that for all $x,y\in \Si$:
\ben
\item $\supp(xy) = \supp(x) \vee \supp(y)$
\item  $\supp(x) \le \supp(y) \quad \Leftrightarrow \quad yx=y $
\een 
\eprop

We call $\La$ the {\em support lattice} and $\supp$ the {\em support map}. Observe that
$$\supp^{-1}(\,\hz\,)=\{e\} \;\;\mbox{ and }\;\; \supp^{-1}(\,\ho\,)= \mSi,$$
where $\hz$ and $\ho$ denote the bottom and top elements of $\La$. In fact, the following
conditions on an element $c\in\Si$ are equivalent:
\ben
\item $\supp(c)=\ho $,
\item  $c\in \max(\Si)$,
\item $cx=c$, for all $x\in \Si$.
\een

Here is the main result on the random walks of Definition \ref{randomwalk}.

\bthm[Brown \cite{Bro1}] \label{Bromain}
Let $\Si$ be an LRB semigroup and $\La$ its support lattice. 
Furthermore, let $\{w_x\}$ be a probability
distribution on $\Si$ and $P_w$ the transition matrix of the induced random walk on the ideal
$\max(\Si)$:
$$P_w(c,d) =\sum_{x\st xc=d} w_x
$$
for $c,d\in \mSi$. Then,
\ben
\item The matrix $P_w$ is diagonalizable.
\vspace{1mm}
\item For each $X\in\La$ there is an eigenvalue
$ \varepsilon_X = \sum_{y\st \mathrm{supp} (y)\le X} w_y \, .
$
\vspace{4mm}
\item The multiplicity of the eigenvalue $\varepsilon_X$ is
$m_X = \sum_{Y\st Y\ge X} \mu_{\La}(X,Y) c_Y ,$\\[2mm]
where $c_Y \defeq {| \max(\Si_{\ge y}})|$,
for any $y\in \supp^{-1}(Y)$.
\vspace{3mm}
\item These are all the eigenvalues of $P$.
\vspace{1mm}
\item Suppose that $\Si$ is generated by $\{x\in \Si\st w_x >   0\}$. Then the random walk on
$\mSi$ has a unique stationary distribution $\pi$.
\een
\ethm

By M\"obius inversion the multiplicities can be determined also from the relations
\beq \label{multipl}
c_X =\sum_{Y\st Y\ge X} m_Y .
\eeq

Theorem \ref{Bromain} is a generalization from the special case of face semigroups
of real hyperplane arangements, to be briefly reviewed in the following section.
In that case the theorem emanates from the work of Bidigare,
Hanlon and Rockmore \cite{BHR} and was expanded by Brown and Diaconis \cite{BrDi}.
The generalization to LRB semigroups was given by Brown \cite{Bro1, Bro2}.

The cited papers also contain information about the rate of convergence to stationarity
and various descriptions of the stationary distribution, e.g. via sampling techniques,
%These parts of the results will not be touched upon in this paper, we refer to the
see \cite{BHR, Bro1, Bro2, BrDi, Dia} for such information.

The following proposition describes two ways in which smaller 
LRB semigroups are induced.

\begin{prop}[\cite{Bro1}] \label{fiber}
Let $\Si$ be an LRB semigroup with support lattice $\La$. Suppose that $x\in \Si$ and
$X\in \La$. Then
\ben
\item $\Si_{\ge x}\defeq \{y\in \Si\st y\ge x\}$ is an LRB semigroup 
%(we call it the {\em filter semigroup} at $x$), 
whose
support lattice is the interval  $[\supp(x), \ho\,]$ in $\La$.
\item If $\supp(x)=\supp(y)$ then $\Si_{\ge x}\cong \Si_{\ge y}$.
\item $\mathrm{Fib}_{\La}(X)\defeq   \{y\in \Si\st \supp(y)\le X\}$ is an 
LRB semigroup (we call it the {\em fiber semigroup} at $X$), whose
support lattice is the interval  $[\,\hz, X]$ in $\La$.
\een
\end{prop}

\subsection{Walks on $\R$-arrangements} \label{sect4:Rarr}

Let $\AA$ be an essential hyperplane arrangement in $\R^d$
%As in Section \ref{sect2:basics} let  be the
with face semilattice $F_{\AA}$ and intersection lattice $L_{\AA}$. 
The following is easily seen from observations % (\ref{semigroup})
(2.1) -- (2.5).

\bprop\label{LRBsemigroup}
$(F_{\AA}, \circ ) \mbox{ is an LRB semigroup} $ with support lattice $\LA$
and support map $\spa$.
\eprop

Let $C_{\AA}$ be the set of regions in the complement of  $\AA$.
There is a one-to-one correspondence $C_{\AA}\leftrightarrow \max(F_{\AA})$.
Thus the general theory produces  a class of random walks on $C_{\AA}$
to which Theorem \ref{Bromain} is applicable. The description of this case is as follows.

\begin{quote}
{\em Random walk on $C_{\AA}$:
Fix a probability distribution $w$ on $F_{\AA}$.
If the walk is currently in region $C \in C_{\AA}$, then choose a face $X\in F_{\AA}$ according 
to $w$ and move to the region $X\circ C$.}
\end{quote}

Let $P_w$
% = (m_{C, D})_{\,C, D \in C_{\AA}}$ 
be the transition matrix
$$P_w(C, D) = \sum_{F\st F\circ C=D} w_F$$
Theorem \ref{Bromain} specializes to the following, where part (3) 
relies on Zaslavsky's formula (Theorem \ref{zaslavsky})
%$$\mbox{number of  regions $=$ size of matrix }=\sum_{F\in L_{\AA}} |\mu (\hz, F)|$$
together with relation (\ref{multipl}).

\begin{thm}[Bidigare-Hanlon-Rockmore \cite{BHR}, Brown-Diaconis \cite{BrDi}] .
\ben\label{Rbrown}
\item $P_w$ is diagonalizable.
\item For each $X\in \LA$ there is an eigenvalue
$\ep_{X}=\sum_{F \st \spa(F) \subseteq X} w_F. $
\item The multiplicity of $\ep_X$ is $|\mu_{\LA}(X, \ho\, )|$. 
\item These are all the eigenvalues.
\item 
Assume that the probability mass $w$ is not concentrated on any single 
hyperplane $H_i$. 
Then there is a unique stationary distribution $\pi$. 
\een\end{thm}

\brem {\rm 
The following interesting result appears in \cite{BBD}. Let  $w$ be the uniform
distribution on the set of vertices (minimal 
elements of $\FA\setminus\{0\}$) of an arrangement in $\R^3$.
Then the probability (according to $\pi$) of being in a region with $k$ sides is proportional to $k-2$.
It is an open problem to give any such geometric characterization of 
the stationary distribution %$\pi$  for 
for arrangements in $\R^d$, $d\ge 4$. }
\erem

\subsection{Walks on $\C$-arrangements}  \label{sect4:Carr}

Let $\AA$ be an essential hyperplane arrangement in $\C^d$
%As in Section \ref{sect2:basics} let  be the
with face semilattice $F_{\AA}$ and intersection lattices $L_{\AA}$ and $\LAaug$.
The following strengthening of observation (\ref{semigroup2})
is immediate.

\bprop\label{LRBsemigroupC}
$(F_{\AA}, \circ ) \mbox{ is an LRB semigroup} $ with support lattice $\LAaug$
and support map $\spa$.
\eprop
\smallskip

Applying the general  theory directly to $\FA$ and the ideal $\max(F_{\AA})$ we
get a walk on the maximal complex sign vectors which is a direct analogue of the real walks
in Section \ref{sect4:Rarr}.

\begin{thm}
 The statements of Theorem \ref{Rbrown} are 
 valid for the complex walks, with the following replacements
for items (2) and (3):
\ben\label{Cwalkthm}
\item[(2)] For each $X\in \LAaug$ such that $\mu_{\LAaug} (X, \ho\,)\neq 0$
there is an eigenvalue
$\ep_{X}=\sum_{F \st \spa(F) \subseteq X} w_F. $
\item[(3)] The multiplicity of $\ep_X$ is $|\mu_{\LAaug}(X, \ho\, )|$. 
\een\end{thm}

The proof of part (3) relies here on the generalized Zaslavsky formula (Theorem \ref{genzaslavsky})
%$$\mbox{number of  regions $=$ size of matrix }=\sum_{F\in L_{\AA}} |\mu (\hz, F)|$$
together with relation (\ref{multipl}). Note that in the formulation of Theorem \ref{Rbrown}
we need not demand that $\mu(X, \ho\,)\neq 0$, since that is automatically true
for geometric lattices. However, in Theorem \ref{Cwalkthm} all we know is that
the lattice is lower semimodular, which implies that the \mob\ alternates in sign but not
that it is nonzero.

Specializing in various directions there are several semigroup-induced random walks
coming out of this situation. We describe two of them.
\smallskip

{\bf Case 1.} Suppose that the probability mass $w$ is concentrated on the real sign vectors and let
$Z=X\circ iY\in \FA$, for real sign vectors $X$ and $Y$. 
Choose $W\in \FA\cap \{0,+,-\}^t$ according to $w$ and move to
$W\circ Z = (W\circ X)\circ iY$. Then $Z$ and $W\circ Z$ have the same imaginary part $iY$.
It can be checked that the subset of $\FA$ consisting of sign vectors with fixed
imaginary part $iY$ is an LRB semigroup. Note that it doesn't come from a filter
of a fiber, as in Proposition \ref{fiber}.

For complexified real arrangements, where sign vectors correspond to intervals, we have
in this case that
$$ [0,X] \circ [R,S]  = [0\circ R, 0\circ R\circ X\circ S] = [R, R\circ X\circ S] $$
So, probability mass concentrated on elements $[0,X]$ (real sign vectors)
gives a random walk on
the set of intervals $[R,S]$, $S$ maximal, {}{for any fixed element $R$}.
\smallskip

{\bf Case 2.} Let $\AA^{\C}$ be the complexification % in $\C^d$
of a real arrangement $\AA^{\R}$ .
%$\LAaug$
We have that $L_{\AA^\C, \text{}aug} \cong \text{Int}(L_{\AA^{\R}})$. The purpose here is to
determine the transition matrix eigenvalues for the fiber semigroup
$\mathrm{Fib}(X) =   \{y\in F_{\AA^{\C}} \st \supp(y)\le X\}$, for
$X = [\pi, \ho] \in  \text{Int}(L_{\AA^{\R}})$.
The support lattice of Fib$(X)$ is the interval  $[\,\hz, X]$ in
 $L_{\AA^\C, \text{}aug}$, cf. Proposition \ref{fiber}.

%Reversing the partial order, we now  work instead with the interval
%$[X, \ho]=[[\hz, \pi], [\ho, \ho]]$ in $\text{Int}(L_{\AA^{\R}})$. Be aware that
%when reversing the order the roles of $\hz$ and $\ho$ are switched.

Theorem \ref{Bromain} shows that the eigenvalues are indexed by intervals
$[\,\al, \be\, ]\in [[\,\hz, \hz\, ], [\pi, \ho\,]]$, i.e., intervals $[\al, \be]$ such that $\al\le\pi$.
Furthermore, the multiplicity of such eigenvalue is, 
according to Theorems \ref{Cwalkthm},
\ref{genzaslavsky} and \ref{muintervals}, the absolute value of
$$\mu_{\mathrm{Int}(L)} ([\,\al, \be\,], [\pi, \ho\,])=\bca
\mu_L (\al, \pi)\, \mu_L (\be, \ho\,), & \mbox{ if } \pi \le \be\\
0, & \mbox{otherwise.}
\eca
$$
Thus, eigenvalues of positive multiplicity occur only when $\al\le\pi\le\be$,
and we have proved the following.

\begin{thm} \label{fiberwalks}
 The statements of Theorem \ref{Bromain} are 
 valid for the complex hyperplane walks induced on fibers {\rm Fib}$(X)$,
 as explained,
 with the following replacements
for items (2) and (3):
\ben\label{Cwalkthm2}
\item[(2)] For each $(\al, \be) \in [\hz, \pi] \times [\pi, \ho] $
there is an eigenvalue
$\ep_{(\al, \be)}$.  %=\sum_{F \st \spa(F) \subseteq X} w_F. $
\item[(3)] The multiplicity of $\ep_{(\al, \be)}$ is 
$|\mu_L (\al, \pi)\, \mu_L (\be, \ho\,)|$.
\een\end{thm}

The exact value of $\ep_{(\al, \be)}$ can of course be stated
as  a special case of Theorem \ref{Bromain}, but we leave this aside.

%$[\hz, X]=[[\hz, \hz], [\pi, \ho]]$

\subsection{Walks on  libraries}  \label{sect4:library}

This section concerns the walks produced by the braid arrangements, both
real and complex. By translating from permutation and partition structures
%that appear in connection with the braid arrangement
we can interpret the states of such walks as distributions of books on 
shelves. This library terminology also 
provides a convenient image for picturing and explaining these walks.
%More precisely, the shelves are ordered from top to bottom (say), and
%the books on each shelf are ordered from left to right. Call such a state a {\em library
%configuration.}

%Different walks are possible, depending on how one specializes the theory.
 %For instance,
%a step in a Markov chain (driven by some probability distribution)
%can be a reordering of
%the shelves plus a reordering of the books on each shelf according to some
%combinatorial rule.
\smallskip

{\bf Real case.}
Here one obtains random walks on permutations governed by probability distributions $w$ on
ordered partitions. This case is thoroughly discussed and exemplified in the literature,
see \cite{BHR, BBD, Bro1, BrDi, Dia}. We mention just two examples.

First, suppose that the probability mass is concentrated on the two-block ordered partitions
whose first block is a singleton. That is,
$$\mbox{probability }  =
\bca w_i,  \mbox{ for the partition } \{i\}\mid [n]\setminus \{i\} \\
0,  \mbox{ \; for all other ordered partitions.}
\eca $$
Then the random walk is precisely the Tsetlin library, for which book $i$
is chosen with probability $w_i$ and moved to the begining of the shelf.

Second, more generally allow non-zero probability for all two-block ordered partitions:
$$\mbox{probability }  =
\bca w_E,  \mbox{ for the partition }  \, E\mid [n]\setminus E \\
0,  \mbox{ \; for all other ordered partitions.}
\eca $$
Then the steps of the random walk consist of removing the books 
belonging to the subset $E$ 
with probability $w_E$
and then replacing them in the induced order at the beginning of the shelf.

In the general case, when non-zero probability is allowed for arbitrary ordered partitions, 
we obtain the one-shelf dynamic library with several borrowers described in the Introduction.

%\begin{thm}[Bidigare-Hanlon-Rockmore, Brown-Diaconis].\\
%Assume $w$ is not concentrated on any single $H_i$. \\
%(a) There is a unique stationary distribution $\pi$. \\
%(b) Algorithm how to sample a region distributed from $\pi$.\\
%(c) Measure of convergence  to $\pi$ of Markov chain.
%\end{thm}

\smallskip

{\bf Complex case.} Let us now see what happens in the case of the complex braid arrangement.
We work out the case of a particular fiber LRB, namely the one
determined by choosing $X=[\pi, \ho\,]$, where $\pi$ is a
partition $(B_1, \ldots, B_k)\in \Pi_n$
and $\ho\,$ is the partition into singletons. 

In our library there are $n$ books labeled by the integers $1$ through $n$,
and $k$ shelves labeled by the integers $1$ through $k$.
Think of $\pi$ as a division of the books into $k$ groups corresponding to the
blocks $B_i$. For instance, $B_1$ could be the set of books on combinatorics, $B_2$
the set of algebra books, and so on. We are going to consider
placements of these $n$ books on the $k$ shelves so that the books in any
particular class $B_i$ stand (in some order) on some particular shelf
dedicated to that class. 

The inverse image $\supp^{-1}(X)$ consists of pairs $[p,s]$, where $p$ is an
ordered partition of the given blocks, $p= \lang   B_{p_1}, \ldots, B_{p_k}\rang $, 
and $s$ is a permutation of $[n]$ refining $p$.
We interpret such an element $[p,s]$ as a particular placement of
the  books: the books in $B_{p_1}$ stand on the top shelf
in the order assigned by $s$, then the books in $B_{p_2}$ stand on the next shelf
in the order assigned by $s$, and so on.

The fiber semigroup Fib$(X)=\supp^{-1}(\La_{\le X})$ consists of pairs $[q,t]$, where $q$ is an
ordered partition such that $\supp(q)$ is a coarsening
of the given partition $\pi= \{B_1, \ldots, B_k\}$, 
and $t$ is an ordered partition refining $q$.

A step in the Markov chain is of the form
$[p,s] \mapsto [q,t] \circ [p,s] = [q\circ p, q\circ p\circ t\circ s]$.
What is its combinatorial meaning? Well, $q\circ p$ is an
ordered partition with blocks $B_1, \ldots, B_k$, and 
$q\circ p\circ t\circ s$ is a permutation refining $q\circ p$. Hence, the combinatorial meaning
of such a step in the Markov chain is that we permute the shelf assignments for the
blocks $B_i$ according to $q\circ p$, and then permute the books on each shelf
as induced by the permutation $q\circ p\circ t\circ s$.

Here is a concrete example. Say we have $14$ books of $4$ types, namely the
algebra books $B_{\mbox{alg}}=\{1,4,5,7\}$, the combinatorics books
$B_{\mbox{comb}}=\{2,8,11,12,14\}$, the geometry books
$B_{\mbox{geom}}=\{6,13\}$, and the topology books
$B_{\mbox{top}}=\{3,9,10\}$. Furthermore, say that the present state of the Markov
chain is this library configuration:
\beq\label{state1}
\barr{ccccc}
11 &  14 &  2 &   12  & 8\\
\hline
6 &   13&&&\\
\hline
4&   7&   5  & 1&\\
\hline
10 &  9 &  3&&\\
\hline
\earr
\eeq
So, in particular,  we have the combinatorics books on the top shelf, the geometry books on the next shelf,
and so on.

Now, let $$q=\lang B_{\mbox{alg}} \mid B_{\mbox{comb}} \cup 
B_{\mbox{top}} \mid  B_{\mbox{geom}}\rang  $$ and
$$t= \lang  4,5  \mid 1,7 \mid8,9,12\mid 14 \mid2,3,10,11 \mid6,13 \rang $$
Then, $[q,t]$ acting on the state (\ref{state1}) leads to  the following configuration
\beq\label{state2}
\barr{ccccc}
4&   5&   7  & 1&\\
\hline
12 &  8 &  14 &  11  & 2\\
\hline
9 &  10 &  3&&\\
\hline
6 &   13&&&\\
\hline
\earr
\eeq

From now on we specialize the discussion to what seems like  a ``realistic'' special case,
in which the Markov chain is driven by choices of subsets  $E\subseteq[n]$ of the books.
This walk has the following description in words.

\begin{quote}{\em Library walk:
A borrower enters the library
and borrows a subset $E\subseteq[n]$ of the books with probability $w_E$. 
These books may come from several shelves.
When returned the books are put back
in the following way. Permute the shelves so that 
the ones that contained one of the borrowed books %from $E$
become the top ones,  maintaining the induced order among them and among the
 remaining shelves, which are now at the bottom. Then, on each shelf place the books 
 belonging to $E$ at the beginning of the shelf,  in the
 induced order, followed by the remaining books in their induced order.
}\end{quote}

For example, 
if this procedure is carried out on the library configuration (\ref{state1})
for the choice $E=\{1,2,3,4\}$ we obtain the new configuration (\ref{state3}).

\beq\label{state3}
\barr{ccccc}
2 &  11 &  14 &  12  & 8\\
\hline
4&   1&   7  & 5&\\
\hline
3 &  10 &  9&&\\
\hline
6 &   13&&&\\
\hline
\earr
\eeq
\smallskip

%Let $K_E$ be the collection of all books on the shelves that are disturbed. Turning
In mathematical language, the following is going on.
For the subset $E\subseteq[n]$ let  $K_E  \defeq \cup_{i\st B_i \cap E \neq \emptyset} \,B_i$ and 
$$q_E\defeq \lang K_E\mid \,[n]\setminus K_E\rang  
\mbox{ and }  t_E\defeq\lang   E\mid K_E\setminus E\mid [n]\setminus K_E \rang .$$
%Now put the probability mass on only elements of this type.
The mathematical description of the library walk is that
we assign the following distribution 
$$\mbox{probability }  =
\bca w_E,  \mbox{ for the partition interval } [q_E, t_E], \text { all }E\subseteq [n] \\
0,  \mbox{ \; for all other intervals of ordered partitions.}
\eca $$
to the elements of the fiber semigroup Fib$([\pi, \ho\,])$, and then we refer to
Theorem \ref{fiberwalks} for the consequences.

To exemplify how the interval $[q_E, t_E]$  acts on a library configuration we 
return once more  to the configuration (\ref{state1}).
Suppose that $E=\{1,2,3,4\}$ and let the interval
$[q_E, t_E]$  act on (\ref{state1}). This leads to the library configuration
(\ref{state3}),  as shown by the following computation table

\beqs\label{state6}
\barr{c||c|c|c|c|c||c|c||c|c|c|c||c|c|c||}
\circ &11  &14  &2&12 & 8 & 6 & 13 & 4 & 7 &5  &1  &10  &9 & 3\\
\hline\hline
1,2,3,4&  &  &2& &  &  &  & 4 &  &  &1  &  & & 3\\
\hline
5,7,8,9,10,11,12,14&11  &14  &&12 & 8 &  &  &  & 7 &5  &  &10  &9 & \\
\hline\hline
6,13 &  &  && &  & 6 & 13 &  &  &  &  &  & & \\
\hline\hline
\earr
\eeqs
\smallskip

Summing up the discussion we obtain the following result.

\begin{thm} \label{fiberwalks2}
 The statements of Theorem \ref{fiberwalks} are 
 valid for the library walk,
 with the following replacements
for parts (2) and (3):
\ben\label{Cwalkthm3}
\item[(2)] For each pair of unordered partitions $(\al, \be)$% \in [\hz, \pi] \times [\pi, \ho] $
such that $\al\le\pi\le\be$ %in $\Pi_n$
(i.e., $\be$ refines $\pi$ and $\pi$ refines $\al$)
there is an eigenvalue
$\ep_{(\al, \be)}$. Furthermore,
$$\ep_{(\al, \be)}=\sum w_E ,
$$ 
the sum extending over all $E\subseteq [n]$ such that $E$ is a union of blocks from $\be$ 
and the shelves containing some element of $E$ is a union of blocks from $\al$.
 %=\sum_{F \st \spa(F) \subseteq X} w_F. $
\item[(3)] The multiplicity of $\ep_{(\al, \be)}$ is 
$\prod (p_i-1)! \prod (q_j-1)! $, where $(p_1, p_2, \ldots)$
are the block sizes of $\be$ and $(q_1, q_2, \ldots)$ the block sizes of $\al$ modulo $\pi$.
\een\end{thm}

Here part (3) uses the well-known formula for the \mob\ of the partition lattice $\Pi_n$
in terms of factorials, see e.g. \cite[p. 128]{EC1}

\bex{\rm
We exemplify the preceding with a worked-out example. Let $n=3$ and 
$\pi = (1,2 \mid 3).$ Then there are four library configurations indexing
the rows and columns of the transition matrix $P_w$:
\newcommand{\confa}{\begin{tabular}{cc}
%\vspace{.01mm}
1&2\\
\hline
3 & \\
\cline{1-1}
\vspace{.01mm}
\end{tabular}}
\newcommand{\confb}{\btab{cc}
2&1\\
\hline
3 &\\
\cline{1-1}
\vspace{.01mm}
\etab}
\newcommand{\confc}{\btab{cc}
3&\\
\cline{1-1}
1&2 \\
\hline
\vspace{.01mm}
\etab}
\newcommand{\confd}{\btab{cc}
3\\
\cline{1-1}
2&1\\
\hline
\vspace{.01mm}
\etab}

\bce
$
\barr{l|c|c|c|c|}
&\confa&\confb&\confc&\confd \\
\hline
\confa &w_1+w_{1,2}+w_{1,3}&w_1+w_{1,3}&w_1+w_{1,2}&w_1    \\
\hline
\confb &w_2+w_{2,3}&w_2+w_{1,2}+w_{2,3}&w_2&w_2+w_{1,2}    \\
\hline
\confc &w_3&0& \quad w_3+w_{1,3} \quad &w_{1,3}     \\
\hline
\confd &0&w_3&w_{2,3}& \quad w_3+  w_{2,3} \quad  \\
\hline
\earr
$
\ece
We ignore the trivial choices $E=\emptyset$
and $E=\{1,2,3\}$, so six elementary probabilities $w_E$ are
assigned. For instance, the entry $w_2 + w_{1,2}$ records that
if books $E$ are removed from the library configuration \;\raisebox{-1mm}{\confd}\; 
and replaced according to the rules, then configuration 
\;\raisebox{-1mm}{\confb}\; 
is obtained precisely if $E= \{2\}$ or $E= \{1,2\}$.

We have that $\hz \cov \pi \cov \ho$ ($\cdot\cov\cdot$ indicates coverings),
so according to Theorem \ref{fiberwalks2} 
there are four pairs $(\al, \be)$ indexing the eigenvalues, all of which have
multiplicity one, and these eigenvalues are
$\bca \ep_{(\,\hz, \pi)}=0 \\
\ep_{(\,\hz, \ho\, )}=w_{1,3}+w_{2,3} \\
\ep_{(\,\pi, \pi\, )}=w_{3}+w_{1,2} \\
\ep_{(\,\pi, \ho\, )}=1
\eca
$

\vspace{3mm}
\noindent
It is instructive to also check how
the elementary probabilities $w_E$ contribute to the various eigenvalues
$\ep_{(\al, \be)}$ in terms of the associated intervals:
%\vspace{6mm}
\bce
\btab{|c|c|c|}
$E$ & $[ q_E, t_E ]$ & contributes to $\ep_{(\al, \be)}$ \\
\hline
$1$ & $[ \lang 12 \mid 3 \rang , \lang 1\mid 2\mid 3 \rang ]$ & $[\al, \be]=[\pi, \ho\,]$ \\
$2$ & $[ \lang 12 \mid 3 \rang , \lang 2\mid 1\mid 3 \rang ]$ & $[\al, \be]=[\pi, \ho\,]$ \\
$3$ & $[ \lang 3 \mid 12 \rang , \lang 3\mid 12 \rang ]$ & $[\al, \be]=[\pi, \ho\,] \;\text{or}\; [\pi, \pi]$ \\
$1,2$ & $[ \lang 12 \mid 3 \rang , \lang 12 \mid 3 \rang ]$ & $[\al, \be]=[\pi, \ho\,] \;\text{or}\;  [\pi, \pi]$ \\
$1,3$ & $[ \lang 123 \rang , \lang 13\mid 2 \rang ]$ & $[\al, \be]=[\pi, \ho\,] \;\text{or}\; [\hz, \ho]$ \\
$2,3$ & $[ \lang 123 \rang , \lang 23\mid 1 \rang ]$ & $[\al, \be]=[\pi, \ho\,] \;\text{or}\; [\hz, \ho]$ 
\etab
\ece

}\eex

\subsection{Walks on  greedoids}  \label{sect4:greedoids}

Denote by $E^*$ the set of repetition-free words $\al=x_1 x_2 \ldots x_k$ in
letters $x_i\in E$. A {\em greedoid} is a language $\LL\subseteq E^*$ such that

\ben
\item[(G1)] if $\al\be\in\LL$ then $\al\in\LL$, for all $\al, \be\in E^*$,
\item[(G2)] if $\al, \be\in \LL$ and $|\al |>  |\be|$, then $\al$ contains a letter $x$ such that
$\be x\in \LL$.
\een
The words in $\LL$ are called {\em feasible} and the longest feasible words 
are called {\em basic}. All basic words have the same length, and 
 $\LL$ is  determined by the basic words as the collection of all their prefixes.

Greedoids were introduced in the early 1980s by Korte and Lov\'asz, see
the accounts in \cite{BjZi2} and \cite{KLS}. The concept can equivalently be formulated 
in terms of set systems, but only the (ordered) language version will concern us here.

Important examples of greedoids are provided by matroids (abstraction of linear hull) and antimatroids 
(abstraction of convex hull). Other examples come from branchings in rooted
directed graphs %(this will be described below) 
and various optimization procedures
(involving some versions of ``the greedy algorithm'').

If $\al, \be\in \LL$ and $|\al |>  |\be|$, then repeated use of 
the exchange property (G2) shows that
$\be$ can be augmented to a word $\be x_1 x_2 \ldots x_j$ with $j=
|\al |-|\be|$ letters $x_i$ drawn from $\al$. But the letters $x_i$ might not
occur in $\be x_1 x_2 \ldots x_j$ in the ``right'' order, i.e., in the order induced by their placement
in $\al$. This motivates defining an important subclass of greedoids.

\begin{ddef} An {\em interval greedoid} is a language $\LL\subseteq E^*$ satisfying
(G1) and the following strong exchange property
\ben
\item[(G3)] if $\al, \be\in \LL$ and $|\al |>  |\be|$, then $\al$ contains a subword $\ga$
of length $|\ga| =|\al |-|\be|$  such that $\al \ga\in \LL$.
\een
\end{ddef}
By {\em subword} we mean what can be obtained by erasing some letters and then
closing the gaps. Matroids, antimatroids and branchings are examples of
interval greedoids.

Let $\LL$ be a greedoid on the finite alphabet $E$. We define an equivalence relation
on $\LL$ by
\beq
\al \sim \be \quad\Lrarr\quad \{\ga\in E^*\st \al\ga\in \LL\} = \{\ga\in E^*\st \be\ga\in \LL\}.
\eeq
So, $\al$ and $\be$ are equivalent if and only if they have the same set of
feasible continuations. The equivalence classes \, $[\al]\in\LL/\sim$\; are the {\em flats}
of the greedoid, and the {\em poset of flats}
$$\Phi \defeq (\LL/\sim ,\; \le)
$$
consists of these classes ordered by
$$[\al]\le [\be] \quad\Lrarr\quad \al\ga\sim\be, \mbox{ for some } \ga\in E^*.
$$
For instance, the poset of flats of a matroid defined in this way 
is easily seen to be isomorphic to the usual geometric
 ``lattice of flats'' of matroid theory.

The feasible words of a greedoid %(i.e., the words belonging to $\LL$) 
can be composed
in the following manner. If $x_1 x_2\ldots x_j \in \LL$ and $y_1 y_2\ldots y_k\in\LL$ then
%\beq\label{composition}
%x_1 x_2\ldots x_j \circ y_1 y_2\ldots y_k \defeq \mbox{lex}( x_1 x_2\ldots x_j y_1 y_2\ldots y_k),
%\eeq
%where $\mbox{lex}( x_1 x_2\ldots x_j y_1 y_2\ldots y_k)= x_1 x_2\ldots x_j y_{i_1} y_{i_2}
%\ldots y_{i_e}$, $1\le i_1 <    i_2 <    \cdots <    i_e \le k$ 
\beq\label{composition}
x_1 x_2\ldots x_j \,\circ\, y_1 y_2\ldots y_k \defeq x_1 x_2\ldots x_j y_{i_1} y_{i_2}
\ldots y_{i_e} \eeq
where $i_1<   i_2 <   \ldots <    i_e$ is 
the lexicographically first non-extendable increasing sequence such that
$x_1 x_2\ldots x_j y_{i_1} y_{i_2} \ldots y_{i_e} \in \LL$.
Letting $\al=x_1 x_2\ldots x_j $
it is equivalent to say that 
$\al\circ y_1 y_2\ldots y_k =\al y_{i_1} y_{i_2}\ldots y_{i_e}$ 
%with $i_1<   i_2 <   \ldots <    i_k$  is
%with  $1\le i_1 <    i_2 <    \cdots <    i_e \le k$, where the right-hand side is 
%the lexicographically first non-extendable feasible
%subword of $x_1 x_2\ldots x_j y_1 y_2\ldots y_k$, i.e., 
is the word obtained, starting from $\al$, by processing
the letters $y_{i}$ of $ y_1 y_2\ldots y_k$ from left to right and 
adding at the end of the word being formed only those letters $y_i$
whose inclusion preserves feasibility.
%omitting 
%from the word being formed
%any letter $y_{i_d}$ whose inclusion in $\al y_{i_1} y_{i_2} \ldots y_{i_d}$ 
%at that stage would violate membership in $\LL$.
%Note that this will leave all the letters $x_i$ but will typically lead to the omission
%of some of the letters $y_i$.

For instance, consider the greedoid on $E=\{x,y,z,w\}$ whose 14 basic words are the words 
in $E^*$ of length $3$ that do not begin with a permutation of $\{x,y,z\}$ or $\{z,w\}$.
This greedoid is discussed on pp. 290--291 of \cite{BjZi2}. Here are two sample 
computations: $$x \circ yzw=xyw \; \;\;\;\mbox{ and }\;\;\;\;  
(x\circ z)\circ w =xzw\neq xz=x\circ (z\circ w) $$ This example shows that the
composition (\ref{composition}) is not associative, and hence does not in general produce  a semigroup. 
For this reason we must limit the discussion to a smaller class of greedoids.

\bthm\label{greedoidLRB}
Let $\LL$ be an interval greedoid.
%for which the composition (\ref{composition}) 
Then $\LL$ with the composition (\ref{composition})  is an LRB semigroup.
Its support lattice is the lattice of flats $\Phi$, and its support map 
$\LL \rarr \Phi$ sends a feasible word $\al$ to its class $[\al]$.
%{\rm (ii)} The composition is associative for polymatroid greedoids and directed branching greedoids.
\ethm
%straight-forward, the only part that is a bit tricky is proving associativity, see Appendix xx.

%For the definition of polymatroid greedoids see \cite{BjZi2, KLS}. 
%This class contains matroids, undirected 
%branching greedoids and  poset greedoids (linear extensions of posets).
That matroids give rise to LRB semigroups in this way
 was mentioned by Brown \cite[p.891]{Bro1}. 
In the matroid case the result is quite obvious, whereas  for the
general case some details turn out to be a little more tricky.
The proof  is deferred to Appendix \ref{sectApp:greedoids}.

%\smallskip

Being an LRB semigroup means that Brown's  theory of random walks,
summarized in Section \ref{sect4:semigroups}, is applicable. What can be said
about the eigenvalue distribution when specialized
to greedoid walks? 

There is an eigenvalue $\ep_X$ for each $X\in\Phi$ whose value and multiplicity $m_X$ are
determined according to parts (2) and (3) of Theorem \ref{Bromain}. However, 
as Example \ref{branchingex} shows, for greedoids the
multiplicities do not depend only on the structure of the interval $[X,\ho]$ in $\Phi$,
as was the case in the corresponding situation  for real and complex hyperplane walks. 
%This is shown by the following small example of a
%directed branching greedoid of rank $2$ and cardinality $3$: node set $=\{r,1,2\}$
%and edge set $=\{(r,1),(r,2),(1,2)\}$. Here $\Phi$ is isomorphic to the boolean lattice
%of subsets of $\{1,2\}$, and $m_{\{1\}}=1\neq 0= m_{\{2\}}$.
\smallskip

We now illustrate greedoid walks for the important case of branchings.
Let $G$ be a directed rooted graph with node set  $\{r,1,2, \ldots, n\}$ and 
edge set $E$. A {\em branching} is a tree directed away from the root node $r$.
 A subset $R\subseteq \{1,2, \ldots, n\}$
is {\em reachable} if it is the set of nodes of some branching.

The {\em branching greedoid} $\LL_G$ consists of ordered strings of edges such that each
initial segment is a branching. It models common search procedures on graphs.
See  \cite{BjZi2} and \cite{KLS} for more information.

The poset of flats of $\LL_G$ is the lattice $\Phi_G$ of reachable sets ordered by
inclusion. This is, in fact, a join-distributive lattice, see the cited references.
The support map sends a branching to the reachable  set of its nodes.

%In the case of branchings a little more can be said about the eigenvalues of
%the transition matrices than for general interval greedoids. 
According to Theorem \ref{Bromain} there is
an eigenvalue $\ep_X$ associated with every reachable node set $X$. Its value is
the sum of the probabilities for the branchings that reach a subset of $X$, and
its multiplicity is given by
$$m_X =\sum_{Y\st Y\ge X} \mu (X,Y) c_Y
$$
Here $c_X$ is the number of ordered edge sequences feasibly  extending
(any branching reaching) $X$ to a maximal branching.

Since $\Phi_G$ is join-distributive its \mob\ takes the simple form
$$\mu (X,Y)=\bca
(-1)^{|Y|-|X|}, & \text{if the interval is Boolean,}\\
0, & \text{otherwise.}
\eca
$$
For each reachable set $X$, let dom$(X)$ denote the superset of all nodes  that
are either in $X$ or else can be reached from $X\cup \{r\}$ along a single 
edge of $G$. It is clear that every set of nodes contained between $X$ and
dom$(X)$ is reachable, and that the {\em domination set} dom$(X)$ is maximal with
this property. Hence, we get the following simplified expression for the
eigenvalue multiplicity at $X$:
\beq m_X =\sum_{X\le Y\le \text{dom}(X)} (-1)^{|Y|-|X|} c_Y
\eeq

\bex{\rm   \label{branchingex}
The rooted directed graph in Figure 4 gives a branching greedoid of rank $3$ with
$9$ basic words: {\em abc, abd, acb, ace, aec, aed, bac, bad, bda}. 
All subsets of $\{1,2,3\}$ except $\{2\}$ are reachable.
\begin{center}
\psfrag{1}{\Huge $1$}
\psfrag{2}{\Huge $2$}
\psfrag{3}{\Huge $3$}
\psfrag{r}{\Huge $r$}
\psfrag{a}{\Huge $a$}
\psfrag{b}{\Huge $b$}
\psfrag{c}{\Huge $c$}
\psfrag{d}{\Huge $d$}
\psfrag{e}{\Huge $e$}
\resizebox{!}{35mm}{\includegraphics{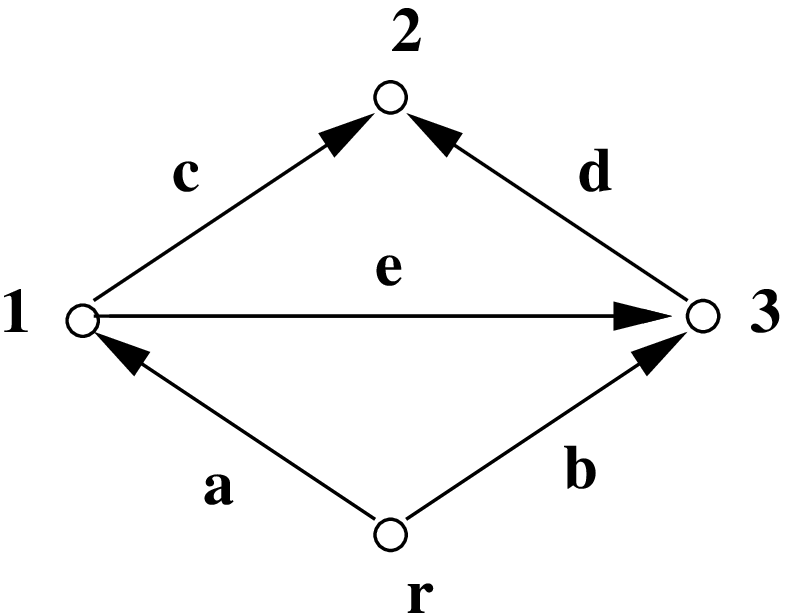}}\\
%{}{Archimedes}
\end{center}
\nopagebreak\vspace{.2cm} \centerline{{\bf Figure 4.} Branching greedoid.}
\vspace{.8cm}
%\noindent
Assign probabilities $w_{\al}$ to the seven feasible words (ordered branchings) of rank one and two:
{\em a, b, ab, ac, ae, ba, bd}.
A step in the random walk on the nine ordered maximal branchings consists in 
choosing one of these words $\al$ according to the
given probabilities $w_{\al}$ and then extending $\al$ to a maximal branching
by adding edges in sequence from the currently visited maximal branching 
according to the composition rule (\ref{composition}).

Here are the eigenvalues for the walk on this  branching greedoid: %in Figure 4:

\begin{center}
$\barr{c||c|c|c|c|}
X & c_X & \text{dom}(X) &m_X & \ep_X \\
\hline
123 & 1 & 123& 1&1 \\
12 & 2 &123& 1 &w_a+w_{ac}\\
13 & 2 &123& 1& w_a+w_b+w_{ab}+w_{ae}+w_{ba}\\
23 & 1 &123& 0& w_b+w_{bd}\\
1 & 6 &123& 3 &w_a\\
3 & 3 &123& 1& w_b\\
\emptyset & 9 &13& 2& 0\\
\earr$
\end{center}
\medskip

}\eex

\brem {\rm 
By copying the procedure that leads from the sign vector system of a real
hyperplane arrangement to that of its complexification
(Section 3.3) we can formally
introduce the complexification of any LRB semigroup. Namely, let $\Si$ be an
LRB semigroup with support lattice $\La$. Define $\Si^{\C}$ to be the set of intervals
$\{[x,y] \st x\le y \text{ in } \Si\}$ with the composition
$$[x,y][z,w] \defeq [xz, xzyw]
$$
One readily verifies that this is an LRB semigroup and that its support lattice
is Int$(\La)$, with the partial order defined in Appendix \ref{sectApp:intervals}.

This way one can complexify e.g. the greedoids walks.
}\erem

\section{Appendix} 

In this section we gather some proofs.
Familiarity with the \mob\ is assumed, a good reference is \cite{EC1}.
%In addition to what can be found there we need for this paper the following results.

\subsection{A generalized Zaslavsky formula} \label{sectApp:zaslavsky}
A ranked poset R with $\hz$ and $\ho$ is said to be {\em Eulerian} if
 $\mu_{R} (x,y)=(-1)^{\rk(y)-\rk(x)}$\; for all $x<   y$ in  $R$. Denote by
 $\max(P)$ the set of maximal elements of a poset $P$.

\bthm\label{genzaslavsky}
Suppose that $f: P\rarr Q$  satisfies the following conditions:
%is an order-preserving map and that the two posets $P$ and $Q$ satisfy
\ben 
\item the posets $P$ and $Q$ are ranked and of the same length $r$,
\item$Q$ has a unique maximal element $\ho_Q$,
\item $\wh{P}\defeq P\uplus\{\hz_P, \ho_P\}$ is Eulerian,
\item $f$ is an order-preserving, rank-preserving and surjective map,
\item[(5)] $\mu_P (f^{-1}(Q_{\le x}))=(-1)^{\rk(x)}$,  for all $x\in Q$,
\item[(6)] $(-1)^{r-\rk(x)}\mu_Q (x, \ho_Q)\ge 0$,  for all $x\in Q$.
\een

Then,
$$  | \max(P) | =\sum_{x\in Q\uplus \hz} |\mu (x, \ho_Q)|
$$
\ethm
\bproof
According to the ``M\"obius-theoretic Alexander duality'' formula
\cite[p. 137]{EC1} condition (3) implies that
$$\mu(R) =(-1)^{r-1} \,  \mu(P\setminus R)
$$
for all subsets $R\sbseq P$. In particular,
\beq\label{alexander}
  | \max(P) |= \mu (\max(P))+1= (-1)^{r-1} \,  \mu(P\setminus \max(P))+1.
\eeq
On the other hand, according to the  ``M\"obius-theoretic fiber formula''  \cite[p. 377]{Wal}
applied to the map $f: P\setminus \max(P) \rarr Q\setminus \ho$ we have that
\beq\label{mufiber} \mu(P\setminus \max(P))  =  \mu(Q\setminus \ho)
- \sum_{x\in Q\setminus \ho} \mu(f^{-1}(Q_{\le x})) \mu (x, \ho_Q) .
\eeq
%The result is obtained by combining (\ref{alexander}) and (\ref{mufiber}) with 
%the input from conditions (5) and (6):
Thus,
\beqarrs
|\max(P)| &=& 1+ (-1)^{r-1} [\, \mu(Q\setminus \ho)
- \sum_{x\in Q\setminus \ho} \mu(f^{-1}(Q_{\le x})) \mu (x, \ho_Q) \,] \\
&=& (-1)^{r} \sum_{x\in Q\uplus \hz} (-1)^{\rk(x)} \mu (x, \ho_Q)) 
= \sum_{x\in Q\uplus \hz} | \mu (x, \ho_Q) | .
\eeqarrs

\eproof

Applying this result to the span map $ \FA \rarr \LA$ of a real hyperplane
arrangement $\AA$ we obtain Zaslavsky's theorem \ref{zaslavsky}.
Applying it to the span map $ \FA \rarr \LAaug$ of a complex hyperplane
arrangement $\AA$ we obtain Theorem \ref{genzaslavskythm}.

%N.B. (1)-(4) implies
%$$  | \max(P) |=(-1)^r \sum_{x\in Q\uplus \hz} 
%\mu_P(f^{-1}(Q_{\le x})\, \mu_Q (x, \ho_Q)
%$$

\subsection{Lattice of intervals}  \label{sectApp:intervals}
Let $L$ be a lattice and Int$(L)\defeq \{(x,y) \st x\le y\}$
the set of its intervals % in $L$
 partially ordered by
$$(x,y) \le (x', y') \quad \mbox{ if and only if \quad $x\le x'$ and }y\le y' .
$$
The poset Int($L$) is itself a lattice with  componentwise operations
$$(x,y) \vee (x', y')= (x\vee x' , y\vee y')  \;\;\mbox{ and }\;\; (x,y) \wedge (x', y')= (x\wedge x' , y\wedge y') .
$$
Its \mob\ is related to that of $L$ in the following way.

\bthm\label{muintervals}
$$\mu_{\mathrm{Int}(L)} ((x,y),(x',y'))=\bca
\mu_L (x,x')\, \mu_L (y,y'), & \mbox{ if } x'\le y\\
0, & \mbox{otherwise.}
\eca
$$
\ethm
\bproof
If $x' \le y$ then $[(x,y),(x',y')] \cong [x,x']\times [y,y']$, so this case follows from the product
property of the \mob.

Assume that $x' \not\le y$. We claim that the element $[x'\wedge y, x'\vee y]$ lacks a lattice-theoretic 
complement in the interval $[(x, y), (x', y')]$.
For, say that $[s,t]$ is such a complement. This means that
$$
\barr{ccccc}
s\vee (x'\wedge y)=x' & &&  &t\vee (x'\vee y)=y'  \\ 
s\wedge (x'\wedge y)=x  &  && &t\wedge (x'\vee y)=y
\earr
$$
Then: \qquad\qquad\qquad $s\le x' \wedge t \le (x' \vee y) \wedge t=y $
\beqarrs
&\Rarr & s \le x' \wedge y \\
&\Rarr & s=s\wedge (x' \wedge y)=x \\
&\Rarr &  x'=x\vee (x' \wedge y)=x' \wedge y \\
&\Rarr &  x' \le y,
\eeqarrs
contradicting the assumption. The interval $[(x, y), (x', y')]$ is not
complemented, so by Crapo's complementation theorem \cite[p. 160]{EC1}
its \mob\ is zero.

\eproof

\subsection{Interval greedoids}  \label{sectApp:greedoids}

The lattice-theoretical structure of semimodularity is closely related to interval greedoids.

\bthm{ \cite[Thm. 8.8.7]{BjZi2}} \label{semimod}
The poset of flats $\Phi$ of an interval greedoid is a semimodular lattice.
Conversely, every finite semimodular lattice arises from some interval greedoid
in this way.
\ethm

This will be used in the proof of Theorem \ref{greedoidLRB},
to which we now turn.
For economy of presentation we assume familiarity with the notation,
conventions and results on pp. 332 -- 334 of \cite{BjZi2}. See particularly the proof of 
Theorem 8.2.5 on p. 334.

\bproof 

\newcommand{\Cov}{<   \!\!\!\cdot\,}

Let $\al =x_1 \ldots x_j$ and $\be =y_1 \ldots y_k$  be feasible words of an interval
greedoid $\LL$. 
By letting  $X_i = [x_1 \ldots x_i]$ and $Y_i = [y_1 \ldots y_i]$,
these words correspond to edge-labeled unrefinable chains
$\emptyset <   \!\!\!\cdot\, X_1 \Cov\: \cdots\: \Cov X_j$ and
$\emptyset <   \!\!\!\cdot\, Y_1 \Cov\: \cdots\: \Cov Y_k$
in the semimodular lattice $\Phi$.
In the same manner (cf. Lemma 8.8.8 of \cite{BjZi2}) the composition 
%$x_1 x_2\ldots x_j \circ y_1 y_2\ldots y_k \defeq \mbox{lex}( x_1 x_2\ldots x_j y_1 y_2\ldots y_k)$
$ x_1 x_2\ldots x_j \circ y_1 y_2\ldots y_k$
corresponds to the edge-labeled unrefinable chain
$$\emptyset <   \!\!\!\cdot\, X_1 \Cov\: \cdots\: \Cov X_j \le X_j \vee Y_1 \le \:\cdots\: \le X_j \vee Y_k.$$
Here, due to semimodularity, the relation $X_j \vee Y_i \le X_j \vee Y_{i+1}$
is either a covering $X_j \vee Y_i \Cov X_j \vee Y_{i+1}$ or an equality $X_j \vee Y_i = X_j \vee Y_{i+1}$,
in which case we omit it from the chain.
This shows that
\beq\label{rel1}
[\al\circ\be]= X_j \vee Y_k =[\al] \vee [\be] 
\eeq
which in turn is used to see that
\beq\label{rel2}[\be] \le [\al] \quad \Leftrightarrow \quad[\al]\vee [\be]=[\al]
\quad\Leftrightarrow \quad [\al\circ\be]=[\al]
\quad\Leftrightarrow \quad \al\circ\be=\al 
\eeq

Thus, once associativity of the composition of feasible words has been established
the proof will be complete. The
other identities required of an LRB semigroup are trivially fulfilled, since
feasible words lack repeated letters.
Relations (\ref{rel1}) and (\ref{rel2}) then show, in view of Proposition \ref{Broprop},
that $\Phi$ is indeed the support lattice of $\LL$ as an LRB  semigroup.

To deal with associativity, let $\ga$
% =z_1 \ldots z_l$ 
be a third feasible word. We want to show that \beq\label{assoc}
(\al\circ\be)\circ\ga=\al\circ(\be\circ\ga) \eeq
By definition
$$(\al\circ\be)\circ\ga=\al\be' \ga '  \mbox{ and } \al\circ(\be\circ\ga)=\al\be' \ga ''$$
where $\be'$ is a subword of $\be$ and $\ga'$ and $\ga''$ are subwords of $\ga$.
Thus it remains to convince ourselves that $\ga'=\ga''$. 
A crucial first step is to show that they are of equal length.

Let $\emptyset <   \!\!\!\cdot\, Z_1 \Cov\: \cdots\: \Cov Z_l$ be
the edge-labelled chain in $\Phi$ corresponding to $\ga$. 
Then $(\al\circ\be)\circ\ga$ corresponds to the chain
$$\emptyset <   \!\!\!\cdot\, X_1 \Cov\: \cdots\: \Cov X_j \le X_j \vee Y_1 \le \:\cdots\: \le X_j \vee Y_k
\le (X_j \vee Y_k) \vee Z_1 \le \:\cdots\: \le (X_j \vee Y_k) \vee Z_l$$
and $\al\circ(\be\circ\ga)$ corresponds to 
$$\emptyset <   \!\!\!\cdot\, X_1 \Cov\: \cdots\: \Cov X_j \le X_j \vee Y_1 \le \:\cdots\: \le X_j \vee Y_k
\le X_j \vee (Y_k \vee Z_1) \le \:\cdots\: \le X_j \vee (Y_k \vee Z_l)$$
Due to associativity of the lattice join operation $\cdot\vee\cdot$ these chains are identical,
and by construction the induced edge-labelings yield
the words $\al\be' \ga '$ and $\al\be' \ga ''$. 
Hence, being related to the same
segment of the common chain, $\ga'$ and $\ga''$ are of the same length. 

We now prove (\ref{assoc}) by induction on the length of the word $\ga$.
Suppose that $\ga =t$ is a single letter. Then $\ga' =\ga''$ since the subwords
of $t$ of length $0$ and $1$ are unique. Hence,
$$(\al\circ\be)\circ t=\al\circ(\be\circ t)$$
%We can now prove that $\ga=\ga'$ by induction on the length of these words.
Suppose now that $\ga=\de t$, meaning that the last letter of $\ga$ is $t$.
Using the induction assumption and the length one case we obtain
%The induction assumption is that $(\al\circ\be)\circ\de=\al\circ(\be\circ\de)$, 
%which is certainly true if $\de$ is the empty word.
\beqarrs
(\al\circ\be)\circ\ga&=& ((\al\circ \be)\circ\de))\circ t =(\al\circ (\be\circ\de))\circ t \\
&=& \al\circ ((\be\circ\de)\circ t) = \al\circ (\be\circ(\de\circ t)) \\
&=& \al\circ(\be\circ\ga)\ .
\eeqarrs

%$$[(\al\circ\be)\circ\ga]= [\al\circ\be]\vee [\ga] = [\al] \vee [\be]\vee [\ga] = [\al]\vee [\be\circ\ga]
%= [\al\circ (\be \circ \ga)] $$

\eproof

\end{document}